\documentclass[11pt, oneside]{amsart}
\pdfoutput=1
\usepackage{amsmath, amssymb, amsthm, amsbsy, bbm, amsfonts, enumitem, color}

\usepackage[all]{xy}
\usepackage[english]{babel}

\usepackage[pdfpagelabels]{hyperref}

\theoremstyle{plain}
\newtheorem{thm}{Theorem}[section]
\newtheorem*{thm*}{Theorem}
\newtheorem{prop}[thm]{Proposition}
\newtheorem{lemma}[thm]{Lemma}
\newtheorem{corollary}[thm]{Corollary}
\newtheorem{property}[thm]{Property}

\newtheorem*{PropA3*}{Proposition \ref{Prop:sec sec}}

\theoremstyle{definition}
\newtheorem{definition}[thm]{Definition}

\theoremstyle{remark}
\newtheorem{remark}[thm]{Remark}

\newcommand{\cC}{\mathcal{C}}
\newcommand{\cA}{\mathcal{A}}
\newcommand{\cF}{\mathcal{F}}
\newcommand{\A}{\mathbb{A}}

\newcommand{\R}{\mathbb{R}} 
\newcommand{\Z}{\mathbb{Z}}
\newcommand{\Q}{\mathbb{Q}}
\newcommand{\define}{\mathrel{\mathop:}=}

\newcommand{\MS}{\mathbb{A}} %Model space 
\newcommand{\App}{\mathcal{A}}%system of apartments 
\newcommand{\seg}{\mathrm{seg}} % segment between two points
\newcommand{\binfinity}{\partial_\App} %building at infninity notation

\newcommand{\RS}{\mathrm{R}}%{\Phi}%{\mathrm{R}}%Root system
\newcommand{\sW}{\overline{W}} %spherical Weyl group	
\newcommand{\aW}{W} %affine Weyl group
\newcommand{\WT}{{\aW}_T}%{\sW\rtimes T} %affine Weyl group w.r.t. T
\newcommand{\Cf}{\mathcal{{C}}_{f}} %fundamental  / Sector 
\newcommand{\bound}{\partial}

\newcommand{\cat}{\mathrm{CAT}}

\numberwithin{equation}{thm}

\begin{document}

\hypersetup{pdfauthor={Schwer},pdftitle={On axiomatic definitions of  non-discrete affine buildings}}.
\title{On axiomatic definitions of non-discrete affine buildings}
\author{Curtis D. Bennett and Petra N. Schwer  \newline With an appendix by Koen Struyve}
\address{Curtis D. Bennett, Department of Mathematics, Loyola Marymount University, 1 LMU Drive, Suite 2700 Los Angeles, CA 90245
\newline Petra N. Schwer, Mathematical Institute, University of M\"unster, Einsteinstrasse~62, 48149 M\"unster, Germany
\newline Koen Struyve, Department of Pure Mathematics and Computer Algebra, Ghent University, Krijgslaan 281, B-9000 Ghent, Belgium}
\email{cbennett@lmu.edu, petra.schwer@wwu.de, kstruyve@cage.ugent.be}

\thanks{The second author was financially supported by the SFB 478 ``Geometric structure in mathematics'' at the University of M\"unster and the DFG Project SCHW 1550/2-1.  The third author is supported by  the Fund for Scientific Research -- Flanders (FWO - Vlaanderen).}
\date{ \today }

\begin{abstract}
In this paper we prove equivalence of sets of axioms for non-discrete affine buildings, by providing different types of metric,  exchange and atlas conditions.  
We apply our result to show that the definition of a Euclidean building depends only on the topological equivalence class of the metric on the model space. The sharpness of the axioms dealing with metric conditions is illustrated in an appendix. There it is shown that a space $X$ defined over a model space with metric $d$ is possibly a building only if the induced distance function on $X$ satisfies the triangle inequality. 
\end{abstract}

\maketitle

\section{Introduction}\label{sec_introduction}

One shortcoming of the definition of non-discrete buildings and affine $\Lambda$-buildings is that often the axioms can be very hard to verify.  This led Parreau \cite{Parreau} to find equivalent axioms for Bruhat-Tits buildings.  In our work we undertake to generalize her work to $\Lambda$-buildings, as well as to present new equivalent sets of axioms. 
While some of Parreau's proofs carry over to $\Lambda$-buildings, any argument that uses compactness, connectedness properties of $\R$
%%% CB 8/25 %%%
or properties specific to the Euclidean metric must be reworked as $\Lambda$ and the $\Lambda$-metric need not have these properties.  

A second complication arises in the case that $\Lambda=\R$. Here the metric on $\R^n$ that is used to define the $\Lambda$-building is different from the Euclidean metric (which is used in the definition of simplicial or $\R$-buildings).  So what importance does the choice of metric play?  We prove that so long as the metrics are equivalent and compatible with the Weyl structure of the model space (in the sense made precise in \ref{cond:metric}), there will be no change in the definition of a $\Lambda$-building.  Moreover, in this case, the induced metric will necessarily satisfy the triangle inequality.

Euclidean buildings, also referred to as non-discrete affine or $\R$-buildings, form one of the prime examples of $\cat(0)$-spaces and were defined by Bruhat and Tits in  \cite{BruhatTits} and \cite{TitsComo} to study Lie-type groups over local fields, as well as fields with non-discrete valuations.   

The first author introduced the more general class of affine $\Lambda$-buildings in \cite{BennettThesis} and \cite{Bennett}, allowing for groups over Krull-valuated fields, that is fields having a valuation taking its values in a totally ordered abelian group $\Lambda$.  Initially, the two main examples were the case of Lie-type groups over the field $F(x,y)$ of rational functions in two variables, with two possible valuations, namely, one into $\Z\times\Z$ lexicographically ordered, and the second onto the subgroup $\{x+y\pi\mid x,y\in\Z\}$ of $\R$.  Care should be taken in the latter case that while geometrically the building can be embedded in a Euclidean building (by thinking of the valuation as landing in $\R$), the different topology of the underlying group can wreak havoc on many of the properties of and standard proofs on Euclidean buildings. Thus the $\Lambda$-building definition generalizes both the notion of non-discrete  $\R$-buildings and that of $\Lambda$-trees.

All the combinatorial information about a simplicial tree is carried by an integer-valued metric on its set of vertices. Generalizing this view a $\Lambda$--tree, as introduced in \cite{MorganShalen}, is a special kind of $\Lambda$-metric space, namely (roughly) the geodesic 0-hyperbolic ones, where the subclass of $\Z$-trees corresponds precisely to the simplicial trees. 
A lot of information about a group may be derived from an action on a $\Lambda$-tree. Compare Chiswell's book \cite{Chiswell} or Morgan's survey \cite{Morgan} for excellent accounts on this topic. 

Kramer and Tent made use of $\Lambda$-buildings in their study of asymptotic cones and their  proof of the Margulis conjecture \cite{KTcones}, \cite{KSTT}. Recently in  \cite{BaseChange} $\Lambda$-buildings have been shown to be functorial in the underlying field. Functoriality easily implies that asymptotic cones of $\R$-buildings are again $\R$-buildings, which was shown with completely different methods by Kleiner and Leeb \cite{KleinerLeeb}. 

\subsection*{The present paper is organized as follows:}

In Section \ref{Sec_Def} we define affine $\Lambda$-buildings and list the properties and axioms in consideration.
After that we will present our main results in Section~\ref{Sec_mr} where we also give an outline of their proofs. 
Detailed proofs are then given in  Sections~\ref{Sec_LA} through \ref{Sec_infinity}.  In \ref{Sec_infinity} we prove the existence of the spherical building at infinity. 
Further explanation concerning the content of these sections is given after the statement of the main result in Section~\ref{Sec_mr}.
In Section~\ref{Sec_WeylCompatibility} we prove that the standard $\Lambda$-metric (or the Minkowski metric) satisfies the Weyl compatibility condition, a connectedness-like condition.
Finally, the Appendix~\ref{Sec_example}  is devoted to the construction of examples of spaces emphasizing the sharpness of axiom (A5). 

\subsection*{Thanks} We are greatly indebted to the anonymous referee for his or her detailed report and valuable comments on the submitted version of this work which lead to a correction and improvement of the main result. 

%%%%%%%%%%%%%%%%%
\section{Definitions and axioms}\label{Sec_Def}

The model apartment of an affine $\Lambda$-building  is defined by means of a totally ordered abelian group $\Lambda$ and a (not necessarily crystallographic) spherical root system $\RS$, where we use the definition of root system as given by Humphreys~\cite{Humphreys}. A root system is called \emph{crystallographic} if all evaluations of co-roots on roots are integers. This is in particular the case if there exists a corresponding affine diagram. The root system $I_2(8)$ is for example not crystallographic but may be used in our setting.

We will now proceed with the definition of the model space of a $\Lambda$-building and continue with defining the $\Lambda$-buildings themselves. We will however not give a detailed introduction to the subject. The reader interested in a more detailed description of the geometric structure and other properties of $\Lambda$-buildings may have a look at Bennett's work \cite{Bennett, BennettThesis}, where these objects were first defined, or at one of the following: \cite{Diss, BaseChange}.

Just as apartments in the geometric realization of Euclidean buildings are isomorphic copies of $\R^n$, the \emph{model space} $\MS$ of an affine $\Lambda$-building can be thought of as a copy of $\Lambda^n$.

We fix a root system $\RS$ and define 
$$
 \MS(\RS,\Lambda) = \mathrm{span}_F(\RS)\otimes_F \Lambda,
$$ 
where $F$ is a sub-field of the real numbers containing all evaluations of co-roots on roots and $\Lambda$ is a totally ordered abelian group with an $F$-module structure. For example, we could take $F=\Q[\{\alpha^\vee(\beta) \,\vert\, \alpha, \beta \text{ roots }\}]$.

The spherical Weyl group $\sW$ associated to $\RS$ acts on $\MS$ by naturally extending its action on $\RS$. 
An \emph{affine Weyl group} $\WT$ acting on $\MS$, is the semi-direct product of $\sW$ by some $\sW$-invariant translation group $T$ of the model space. In the case that $T$ is the entire space $\MS$, we will write $\aW$ instead of $\WT$. 

Elements  in $\WT$ that can be written as $t\circ r_\alpha$ for some $t\in T$ and a reflection $r_\alpha$ in $\sW$ are called  \emph{(affine) reflections} if their fixed point set is non empty. 
The fixed point set $H=H_{t,\alpha}$ of an affine reflection, which splits $\MS$ into two half-spaces,  is called \emph{(affine) hyperplane}.
Here we assume the hyperplane to be part of each of its two half-apartments making them closed subspaces of $\MS$. 
We say that a hyperplane $H_{\alpha,t}$ \emph{separates} two elements $x,y\in \MS$ if $x$ and $y$ are contained in different open halfspaces determined by $H_{\alpha,t}$.  Two hyperplanes are \emph{parallel} if they are translates of one another. In this case they are of the form $H_{\alpha,t}, H_{\alpha,k}$ for some $t,k\in\Lambda$ and $\alpha\in\RS$. It is  shown in ~\cite{Bennett} that hyperplanes are parallel in our sense if and only if they are at bounded distance. 
Associated to a basis $B$ of the root system $\RS$ there is a \emph{fundamental Weyl chamber} $\Cf$. The chamber $\Cf$ is a fundamental domain for the action of $\sW$ on $\MS$. Its images in $\MS$ under the affine Weyl group are \emph{Weyl chambers} (sometimes called \emph{sectors}). If two Weyl chambers $S$ and $T$ contain a common sub-Weyl chamber we call them \emph{parallel} and write $\partial S=\partial T$ for their parallel class. A \emph{Weyl simplex} is a face of a chamber. The smallest face of dimension 0 is called \emph{basepoint}, and a \emph{panel} is a  Weyl simplex of co-dimension one. 

A subset $C$ of $\MS$ is \emph{(closed) $\WT$-convex} (or just \emph{(closed) convex}) if it is the intersection of finitely many half-apartments. The \emph{convex hull}  of a set $Y \subset \MS$ is the intersection of all half-apartments containing $Y$. Weyl simplices, chambers and hyperplanes are all examples of closed convex sets. 

One can endow $\MS$ with a natural $\aW$-invariant metric taking its values in $\Lambda$, and thus making $\MS$ a $\Lambda$-metric space in the following sense: 
A map $d:X\times X \mapsto \Lambda$ on a space $X$ is a \emph{$\Lambda$-metric} if  for all $x,y,z$ in $X$ the following conditions are satisfied
\begin{enumerate}
\item $d(x,y)\geq 0, \forall x,y$, and $d(x,y)=0$ if and only if $x=y$,
\item $d(x,y)=d(y,x)$ and
\item the triangle inequality $d(x,z)+d(z,y)\geq d(x,y)$ holds.
\end{enumerate}
While in principle there exist many different potential $\Lambda$-metrics for $X$, depending on $\Lambda$, the definition of one may be somewhat complicated.  In the case where $\Lambda=\mathbb{R}$, the standard Euclidean metric works, however, in the case where square roots may not exist ($\Lambda=\mathbb{Q}$ or $\mathbb{Z}\times\mathbb{Z}$ for example) things can be more difficult.  One such solution for all $\Lambda$ is to use a modified Minkowski metric (see \cite{Bennett} for details).

We now define generalized non-discrete affine buildings. Throughout the following fix a model space $\MS$ and an affine Weyl group $\WT$.

\begin{definition}\label{Def_LambdaBuilding}
Let $X$ be a set and $\App$, called the \emph{atlas} of $X$,  be a collection of injective charts $f:\MS\hookrightarrow X$.
For each $f$ in $\App$, we call the images $f(\MS)$   \emph{apartments} and define \emph{Weyl chambers, Weyl  simplices, hyperplanes, half-apartments, etc.\ of $X$} to be images of such in $\MS$ under any $f$ in $\App$. The pair $(X,\App)$ is a \emph{(generalized) affine building} (or $\Lambda$-building) if the following conditions are satisfied
\begin{enumerate}[label={(A*)}, leftmargin=*]
\item[(A1)] The atlas is invariant under pre-composition with elements of $\WT$.
\item[(A2)] Given two charts $f,g\in\App$ with $f(\MS)\cap g(\MS)\neq\emptyset$, then  $f^{-1}(g(\MS))$ is a closed convex subset of $\MS$ and there exists $w\in \WT$ with $f\vert_{f^{-1}(g(\MS))} = (g\circ w )\vert_{f^{-1}(g(\MS))}$.
\item[(A3)] For any pair of points  in $X$ there is an apartment containing both.
\end{enumerate}

Given a $\WT$-invariant $\Lambda$-metric $d_{\MS}$ on the model space,  axioms (A1)--(A3) imply the existence of a $\Lambda$-valued distance function on $X$, that is a function $d:X\times X\mapsto \Lambda$ satisfying all conditions of the definition of a $\Lambda$-metric except possibly the triangle inequality.  The distance between points $x,y$ in $X$ is the distance between their pre-images under a chart $f$ of an apartment containing both.

\begin{enumerate}[label={(A*)}, leftmargin=*]
\item[(A4)] Given two  Weyl chambers  in $X$ there exist sub-Weyl chambers of both which are contained in a common apartment. 
\item[(A5)] For any apartment $A$ and all $x\in A$ there exists a \emph{retraction} $r_{A,x}:X\to A$ such that $r_{A,x}$ does not increase distances and $r^{-1}_{A,x}(x)=\{x\}$.
\item[(A6)] Let $f, g$ and $h$ be charts such that the associated apartments pairwise intersect in half-apartments. Then $f(\MS)\cap g(\MS)\cap h(\MS)\neq \emptyset$. 
\end{enumerate}
By (A5) the well-defined distance function $d$ on $X$ satisfies the triangle inequality.
\end{definition}

One problem raised by the definition is that condition (A5) relies on the choice of a $\Lambda$-metric $d_{\MS}$.  However, as alluded to earlier, there are potentially many possible metrics to choose from. To keep our notation manageable, we will assume the metric on $\MS$ is given as part of $\MS$.

The main goal of the present paper is to prove equivalence of certain sets of axioms. Let us therefore collect all properties which are necessary to state the main result. These properties break into three categories: metric conditions, exchange conditions, and atlas conditions.  

We begin with the \emph{metric conditions}.  Axiom (A5) is one such condition since it implies that the global distance function is a metric on the space $X$.  In \cite{Bennett}, this condition is used to prove the existence of a building at infinity $\partial_\App X$  whose simplices are parallel classes of Weyl-simplices in $X$ as defined in \ref{Def:parallel}.  However, in \cite{Parreau} (for $\Lambda=\mathbb{R}$) a proof of the existence of the building at infinity was given that did not require the full power of condition (A5) but only the triangle inequality, (TI) and properties of the Euclidean distance function. We generalize her proof to the general case in Section~\ref{Sec_infinity}, 
%%%  CB 8/25
although doing so requires identifying the special conditions on the $\Lambda$-metric $d$ inherent in the proof in \cite{Bennett}. 
Consequently, there is benefit to having the weaker condition necessary, namely that the induced distance function on $X$ is a metric.
\begin{itemize}
\item[(TI)] (Triangle inequality) The distance function $d$ on $X$, which exists assuming (A1)--(A3), is a metric,  i.e.\ satisfies the triangle inequality.
\end{itemize}

Alternatively one could directly assume that the space under consideration has a spherical building at infinity. 

\begin{itemize}
\item[(BI)] (Building at infinity) The set $\bound X$ of parallel classes of Weyl simplices is a
spherical building with apartments the boundaries $\bound A$ of apartments $A$ in $X$.
\end{itemize}

We next move to the \emph{exchange type conditions}.  One difficulty with generalizing the definition of an $\mathbb{R}$-building to the $\Lambda$-building case was that the totally ordered group $\Lambda$ might not be topologically connected.  In the case of trees (the lowest dimension affine buildings), the move from $\mathbb{R}$-trees to $\Lambda$-trees required the introduction of a $Y$-condition, which is essentially a condition that says that when two paths diverge, the symmetric difference of those paths together with the point of divergence itself forms a path.  For the higher dimensional $\Lambda$-building context, this condition was encapsulated in (A6).  For our purposes, we have two other exchange conditions, each of which is slightly stronger than (A6).  The first is most naturally an exchange condition, and hence we give it that name.

\begin{itemize}
\item[(EC)] (Exchange condition)  Given two charts $f_1,f_2\in\App$ such that $f_1(\MS)\cap f_2(\MS)$ is a half apartment, then there exists a chart $f_3\in\App$ such that $f_3(\MS)\cap f_j(\MS)$ is a half apartment for $j=1,2$.  Moreover, $f_3(\MS)$ is the symmetric difference of $f_1(\MS)$ and $f_2(\MS)$ together with the boundary wall of $f_1(\MS)\cap f_2(\MS)$.
\end{itemize}

Note that the exchange condition can be restated in ``apartment language'' as: Given two apartments $A$ and $B$ intersecting in a half-apartment $M$ with boundary wall $H$, then the set $(A\oplus B)\cup H$ is also an apartment, where $\oplus$ denotes the symmetric difference of $A$ and $B$.\footnote{By definition $A\oplus B =\{x\in A\setminus B \}\cup \{x\in B\setminus A \}$. }

We will also consider the following even stronger exchange condition for which Linus Kramer suggested the name \emph{sundial configuration}. 

\begin{itemize}
\item[(SC)] (Sundial configuration)  Suppose $f_1\in\App$ and $S$ is a Weyl chamber of $(X,\App)$ such that $P=S\cap f_1(\MS)$ is a panel of $S$.  Let $M$ be the wall of $f_1(\MS)$ containing $P$.  Then there exist $f_2 \neq f_3\in\App$ such that $f_1(\MS)\cap f_j(\MS)$ is a half-apartment and $(M\cup S)\subset f_j(\MS)$ (for $j=2,3$). 
\end{itemize}

The sundial configuration can be restated as: Given an apartment $A$ of $X$ and a chamber $c$ in the building at infinity such that $c$ shares a co-dimension one face $p$ (a panel) with $\partial A$ (the boundary of $A$ or the apartment at infinity associated to $A$) but is not contained in $\partial A$, then there exist two apartments $A_1 \neq A_2$ such that $c\in\partial A_i, i=1,2$ and such that $A_i\cap A$ is a half apartment with bounding wall spanned by a panel in $p$. 

The last set of conditions are the \emph{atlas conditions}.  These conditions all state properties of the atlas set $\App$ in terms of containing subsets of $X$.  Thus conditions (A3) and (A4) are atlas conditions.   These atlas conditions typically correspond to statements about objects (two points or Weyl chambers for example) being contained in an apartment (with one exception).  To be more precise, we need some terminology. 

We say that two Weyl simplices $F$ and $G$ \emph{share the same germ} if both are based at the same vertex and if $F\cap G$ is a neighborhood of $x$ in $F$ and in $G$.
It is easy to see that this is an equivalence relation on the set of Weyl simplices based at a given vertex. The equivalence class of an $x$-based Weyl simplex $F$ is denoted by $\Delta_x F$ and is called the \emph{germ of $F$ at $x$}.
The germs of Weyl simplices at a vertex $x$ are partially ordered by inclusion: $\Delta_x F_1$ is contained in $\Delta_xF_2$ if there exist $x$-based representatives $F'_1, F'_2$ contained in a common apartment such that $F_1'$ is a face of $F_2'$. Let $\Delta_xX$ be the set of all germs of Weyl simplices based at $x$.  We note that since the definition of a germ is dependent on the definition of a neighborhood of a point, the notion of germs are necessarily dependent on the (equivalence class of the) metric on $\MS$ with which we start.

A germ $\mu$ of a Weyl chamber $S$ at $x$ is \emph{contained} in a set $Y$ if there exists $\varepsilon\in\Lambda^{+}$ such that $S\cap B_\varepsilon(x)$ is contained in $Y$ where $B_\varepsilon(x)$  denotes the usual $\varepsilon$-ball around $x$.  We are now ready to state the first three of our new atlas conditions.

\begin{itemize}
 \item[(LA)] (Large atlas) Any two germs of Weyl chambers are contained in a common apartment. % (A3')
 \item[(aLA)] (Almost a large atlas) For all points $x$ and all $y$-based Weyl chambers $S$ there exists an apartment containing both $x$ and $\Delta_yS$.
 \item[(GG)] (Locally a large atlas) Any two germs of Weyl chambers based at the same vertex are contained in a common apartment.
\end{itemize}
Note that both (LA) and (aLA) imply (A3).

We say two $x$-based germs are \emph{opposite} if they are contained in a common apartment $A$  and are images of one another under the longest element of the spherical Weyl group (which acts  on the set of (germs of) $x$-based Weyl chambers in $A$).\footnote{It is easy to see that in the case where we know that the germs at $x$ form a spherical building, then opposite germs are also opposite in the usual spherical building sense.} Two Weyl chambers are \emph{opposite at $x$}, if their germs are opposite.  This leads us to our fourth new atlas condition,

\begin{itemize}
 \item[(CO)] (Opposite chambers) Two opposite $x$-based Weyl chambers $S$ and $T$ are contained in a unique common apartment.
\end{itemize}

For our last atlas conditions, we need a metric notion of a ``convex hull'' like object.  Thus, we define the \emph{segment} $\seg_M(x,y)$ of points $x$ and $y$ in a metric space $M$ to be the set of all points $z\in M$ such that $d(x,y)=d(x,z)+d(z,y)$.

\begin{itemize}
 \item[(FCa)] (Finite cover of apartments) For an arbitrary point $z$ in $X$ every apartment $A$  is contained in a finite union of Weyl chambers based at $z$.
 \item[(sFC)] (Strong finite cover) For any pair of points $x$ and $y$, all apartments $A$ containing $x$ and $y$ and an arbitrary germ $\mu$ based at $z$ in $X$,  the segment $\seg_A(x,y)$ is contained in a finite union of Weyl chambers based at $z$ such that each of these Weyl chambers is  contained in a common apartment with $\mu$.
 \item[(sFCa)] (Strong finite cover of apartments) For an arbitrary germ $\mu$ at a point $z$ in $X$, every apartment $A$ is contained in a finite union of Weyl chambers based at $z$ such that
 each of these Weyl chambers is contained in a common apartment with $\mu$.
\end{itemize}

It is clear that (sFCa) implies (sFC) and (FCa). From Lemma~\ref{Lem_sFC-LA} we get that (sFC) implies (LA) and hence (A3). 
Condition (sFC) is the one of actual interest for our main theorem while the other two appear during the proof. We would have liked to replace all these conditions by a weaker finite cover condition which is implied by all the above and reads as follows:
\begin{itemize}
 \item[(FC)] (Finite cover) For all triples of points $x,y$ and $z$ in $X$ and all apartments $A$ containing $x$ and $y$ the segment $\seg_A(x,y)$ is contained in a finite union of Weyl chambers based at $z$.
\end{itemize}
But this would imply that we (at least) had to add (aLA) in equivalence (\ref{a10}) of Theorem~\ref{MainThmB} and to the assumptions of Lemma \ref{Lm:fin cover} and Propositions \ref{Prop_A5} and \ref{Th:tentwo} .

\begin{remark}
Both, the existence of a large atlas (LA) and its local analog (GG)  were introduced by Parreau \cite{Parreau}. Condition (LA) was called (A3') in \cite{Parreau} according to its proximity to axiom (A3) and the abbreviation (GG) probably stood for ``germe - germe''. 
The opposite chamber property (CO) also appeared in \cite{Parreau}, where (CO) stood for ``chambres oppos\'ees''.
The condition (aLA) to almost have a large atlas is ``in between'' (A3) and the existence of a large atlas and suffices for one of the implications in \ref{MainThmB}. 
\end{remark}

%%%%%%%%%%%%
\section{Main results}\label{Sec_mr}

The purpose of this section is to state our main results.

Recall that we say that $(X,\App)$ is \emph{a space modeled on $\MS$} if $X$ is a set together with a collection $\App$ of injective \emph{charts} $f:\MS\hookrightarrow X$ such that $X$ is covered by its charts. That is $X=\bigcup_{f\in\App} f(\MS)$.  Throughout the remainder of the paper we will assume that $(X,\App)$ satisfies conditions (A1)--(A3).

We will now introduce a metric condition we need to impose on the distance function in our main result below.  

Recall that a hyperplane $H$ in $\MS$ separates two elements $x,y\in\MS$ if $x$ and $y$ lie in different half-spaces determined by $H$. That is if $H$ is the  affine hyperplane 
$$
H=H_{\alpha,\lambda}=\{x\in\MS \;\vert\; \frac{(\alpha,\alpha)}{2}\langle x,\alpha^\vee\rangle =\lambda \},
$$ 
then $x$ and $y$ are separated by $H$ if either
$$
\frac{(\alpha,\alpha)}{2}\langle x,\alpha^\vee\rangle < \lambda < \frac{(\alpha,\alpha)}{2}\langle y,\alpha^\vee\rangle
$$
or the same equation holds with the roles of $x$ and $y$ interchanged. Here $(\;,\;)$ stands for the scalar product on the vector space in the definition of the defining root system $\RS$ and $\langle\;,\;\rangle$ for the evaluation of co-roots on roots.    

%%% CB 8/25 %%%
We define a type function on the translation classes of Weyl simplices in $\MS$. As in \cite{Bennett}, we say that the sector panel $\tilde{P}$ of the fundamental Weyl chamber $\Cf$ with $\langle x, \alpha_i^\vee \rangle=0$ for all $x\in \tilde{P}$ is the \emph{fundamental sector panel of type $i$}.  For a general Weyl simplex $\tilde{F}$ in $\Cf$, the \emph{type of $\tilde{F}$} is the union of the types of the panels that contain it. Since any arbitrary Weyl simplex $F$ in $\MS$ is the image of a unique Weyl simplex $\tilde{F}$ of $\Cf$ under some element $w\in \aW$, we define the \emph{type of a Weyl simplex $F$} to be the type of $\tilde F$.  

\begin{definition}\label{cond:metric}
A $\WT$-invariant metric $d$ on $\MS$ is \emph{Weyl compatible} if the following two conditions hold:
\begin{enumerate}
 \item\label{cond:metric1}
 %%% CB 8/25 %%%
  For all Weyl chambers $S$ and all sector faces $P_1$ and $P_2$ of the same type, if there exists $\lambda\in \Lambda$, such that for all sub-faces $P_1' \subset P_1$ there exists a sub-face $P_2'\subset P_2$ such that for $y\in P_2'$ there exists a $x\in P_1'$ with $|d(v,y)-d(v,x)|\le \lambda$ for all $v\in S$, then $P_1$ and $P_2$ are translates of one another.
 \item\label{cond:metric2} Whenever a hyperplane $H$ separates two points $x$ and $y$ in $\MS$, then $\seg_\MS(x,y)\cap H\neq \emptyset$. 
\end{enumerate}
\end{definition}

%%% CB 8/25 %%%
Condition (\ref{cond:metric1}) is the property proved for the standard $\Lambda$-metric in Theorem~2.20 in \cite{Bennett}. This condition guarantees both that in an apartment sector-faces being at bounded distance is equivalent to being translates of one another (see Lemma~\ref{lem:parallel}) and that non-translate sector panels of the same type can be distinguished by their distances from points of any sector. This is used in Theorem~\ref{thm:infty} for the proof of the existence of a spherical building at infinity.  

Condition (\ref{cond:metric2}) is a connectedness type condition,  ensuring  that if a segment is split into a sequence of sub-segments by hyperplanes then each hyperplane contains a point of the segment (see Lemma~\ref{lem:subsegments}).

A metric $d$ on a $\Lambda$-metric space $M$ is \emph{geodesic} if for all $x,y\in M$ there exists an isometric embedding of the interval $[0,d(x,y)]\subset\Lambda$ into $(M, d)$. It is clear that a geodesic real-valued metric satisfies condition~\ref{cond:metric} (\ref{cond:metric2}). Weyl compatibility is easy to see for the Euclidean metric and shown in Section~\ref{Sec_WeylCompatibility} for the standard $\Lambda$-metric introduced in \cite{Bennett}. We leave as an open question whether geodesic $\Lambda$-metrics are necessarily Weyl compatible.

Using explicit constructions and combinatorial properties of links and the building at infinity we prove in Proposition~\ref{Pr:A6 A6'} that (A6) and the exchange condition (EC) are equivalent assuming (A1)--(A4) and (BI).   By Theorem~\ref{thm:infty}, this follows in the case that (A1)-(A4) and (TI) is true (and consequently also when (A5) is true).
By similar arguments we obtain in Proposition~\ref{Pr:A6' A6''} that (EC) and  (SC) are equivalent. Hence we have the following proposition. 

\begin{prop}\label{MainThmA}
Let $(X,\App)$ be a space modeled on $\MS=\MS(\RS,\Lambda)$ such that the metric is Weyl compatible, axioms (A1)--(A4) are satisfied and any one of (A5), (TI), or (BI) is true. Then 
$$
\mathrm{(A6)} \Leftrightarrow \mathrm{(SC)} \Leftrightarrow \mathrm{(EC)}.
$$
\end{prop}

We are ready to state the main theorem. 
\begin{thm}\label{MainThmB}
Suppose $\MS=\MS(\RS,\Lambda)$ is a model space equipped with a Weyl compatible metric and $(X,\App)$ is modeled on $\MS$. Then the following are equivalent: 
\begin{enumerate}
 \item\label{a1} $(X,\App)$ is an affine $\Lambda$-building, that is axioms (A4), (A5) and (A6) are satisfied. 
 \item\label{a2} $(X,\App)$ satisfies (A4), (A5) and (EC). 
 \item\label{a3} $(X,\App)$ satisfies (A4), (A5) and (SC). 
 \item\label{a4} $(X,\App)$ satisfies (A4), (TI) and (A6). 
 \item\label{a5} $(X,\App)$ satisfies (A4), (TI) and (EC).
 \item\label{a6} $(X,\App)$ satisfies (A4), (TI) and (SC).
 \item\label{a7}  $(X,\App)$ satisfies (TI), (GG) and (CO).
 \item\label{a8}  $(X,\App)$ satisfies (GG) and (CO). 
 \item\label{a9}  $(X,\App)$ satisfies (LA) and (CO). 
 \item\label{a10} $(X,\App)$ satisfies (A4), (sFC) and (EC). 
 \item\label{a11} $(X,\App)$ satisfies (A4), (BI) and (A6). 
%\item\label{a12} $(X,\App)$ satisfies (A4), (BI) and (EC). 
%\item\label{a13} $(X,\App)$ satisfies (A4), (BI) and (SC). 
 \end{enumerate}
\end{thm}

By Proposition~\ref{MainThmA} we can replace (A6) with either (SC) or (EC) in item (\ref{a11}) and by the same proposition combined with \ref{Prop_A5} in item (\ref{a10}) condition~(EC) might be replaced by either (A6) or (SC). 

Before discussing the proof, it is worth saying a few words about the relationship of the various conditions.  In particular, statements (\ref{a4})--(\ref{a6}) and (\ref{a11}) correspond to various options on metric and exchange conditions.  Statements (\ref{a7}), (\ref{a8}) and (\ref{a9}) are the stronger versions of the atlas conditions (each having one on germs and one on Weyl chambers) that imply $(X,\App)$ is an affine $\Lambda$-building, allowing us to ignore both metric and exchange conditions.  Finally, statement (\ref{a10}) lets us use the weakest of the atlas conditions (in that (A4) is often easier to show than (CO)) at the cost of a finite cover condition.  Also recall that conditions (LA), (aLA) and (sFC) imply axiom  (A3).

It is also worth mentioning the role that Weyl compatibility of the metric plays.  The first condition shows up whenever we need to work with the building at infinity and are given (TI).   The second condition is used in only two places,  that is in proving  $(\ref{a10})$ implies $(\ref{a2})$ and  in proving $(\ref{a6})$ implies $(\ref{a1})$ (where only the first one is necessary to prove). In both cases we need to prove (A5), i.e.\ apply Proposition~\ref{Th:tentwo} or~\ref{Th:main}. So it would be interesting to see whether one can avoid the Weyl compatibility condition (2).% or whether one can directly deduce $(\ref{a1})$ or $(\ref{a2})$ from $(\ref{a8})$.

Propositions \ref{prop_sc} and \ref{Prop_LA} in Section~\ref{Sec_localStructure} indicate that starting with (\ref{a11}) even further equivalences could be shown. 

We will show the following implications:

$$
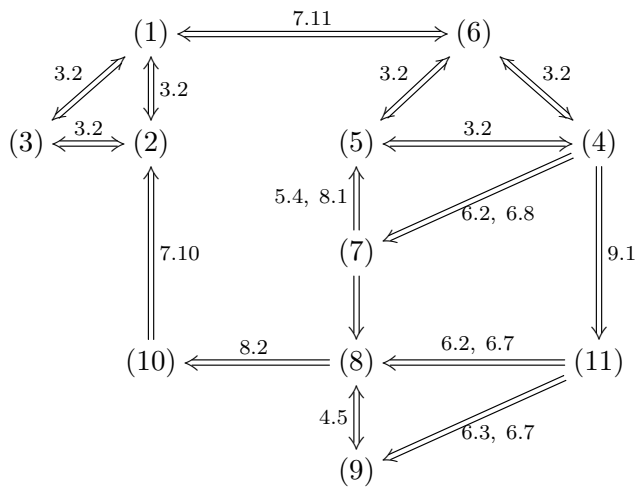
\xymatrix{
& (\ref{a1}) \ar@{<=>}[dl]_{\ref{MainThmA}}\ar@{<=>}[d]^{\ref{MainThmA}}\ar@{<=>}[rrr]^{\ref{Th:main}} &&& (\ref{a6})\ar@{<=>}[dl]_{\ref{MainThmA}}\ar@{<=>}[dr]^{\ref{MainThmA}} \\
(\ref{a3})\ar@{<=>}[r]^{\ref{MainThmA}} & (\ref{a2}) \ar@{<=}[dd]^{\ref{Prop_A5}} && (\ref{a5})\ar@{<=>}[rr]^{\ref{MainThmA}} \ar@{<=}[d]_{\ref{Pr:EC},~\ref{Prop_A4}} && (\ref{a4})\ar@{=>}[dll]^{\ref{Cor_GG},~\ref{Cor_COTI} }\ar@{=>}[dd]^{\ref{thm:infty}} \\
&&& (\ref{a7}) \ar@{=>}[d]\\
&(\ref{a10})\ar@{<=}[rr]^{\ref{eightten}}&& (\ref{a8}) \ar@{<=>}[d]_{\ref{Prop:sec sec}}\ar@{<=}[rr]^{\ref{Cor_GG},~\ref{Cor_CO} }&& (\ref{a11})\\ %\ar@{<=>}[d]_{\ref{MainThmA}}\ar@{<=>}[dr]^{\ref{MainThmA}}\\ 
&&& (\ref{a9})\ar@{<=}[urr]_{ \ref{Prop_LA},~\ref{Cor_CO}}%&& (\ref{a12}) \ar@{<=>}[r]^{\ref{MainThmA}}&(\ref{a13})
}
$$
\medskip

As you can see there is some kind of redundancy in these implications which also clarifies why we included $(\ref{a7})$ in the list of equivalences even though $(\ref{a8})$ is stronger: One can deduce $(\ref{a5})$ from $(\ref{a7})$ which may or may not give room for a possible simplification of the proof of the theorem.   In addition, we have only included Theorem~\ref{thm:infty} where it is explicitly used.  However, it is implicitly used in \ref{MainThmA} as well as in \ref{Thm_residue}, \ref{Cor_COTI} and \ref{Th:main}.

\begin{proof}[Proof of Theorem~\ref{MainThmB}]
By Proposition~\ref{MainThmA} we have that (\ref{a1}), (\ref{a2}) and (\ref{a3}) are equivalent as well as equivalence of (\ref{a4}), (\ref{a5}) and (\ref{a6}). % and (\ref{a11}), (\ref{a12}) and (\ref{a13}). 
Equivalence of (\ref{a1}) and (\ref{a6}) is shown in Proposition~\ref{Th:main}.
That (\ref{a1}) implies  (\ref{a4}) is obvious. 

Assuming (\ref{a4}) we may deduce (BI) from Theorem~\ref{thm:infty}. Therefore (\ref{a4}) implies (\ref{a11}). 
Further we obtain (GG) and  (CO)  as discussed in Corollaries~\ref{Cor_GG}, \ref{Cor_CO} and \ref{Cor_COTI}. Hence (\ref{a4}) implies (\ref{a7}) and item (\ref{a8}) follows from (\ref{a11}).  That (\ref{a7}) implies (\ref{a8}) is clear.

In Section \ref{Sec_LA} we prove Proposition~\ref{Prop:sec sec}, which implies that (\ref{a9}) follows from (\ref{a8}). The converse is obvious.

Theorem~\ref{eightten} proves that (\ref{a8}) implies (\ref{a10}).

Axiom (A5) is verified in Proposition~\ref{Prop_A5} using (A1), (A2) and condition (sFC). Therefore item (\ref{a10}) implies (\ref{a2}).

Finally combining Propositions~\ref{Pr:EC} and \ref{Prop_A4} we obtain that item (\ref{a7}) implies (\ref{a5}) which completes the proof of our main result.
\end{proof}

Part of the power of statements~(\ref{a8}) through~(\ref{a11}) is that they avoid any explicit mention of the distance function on the whole space, which often makes them easier to use. However, this begs the question of why we might heuristically expect them to be equivalent.  It seems that a key issue is condition~(CO) in that by requiring that sectors whose germs are opposite are in an apartment (which is automatically unique by axiom (A2)). Together with (GG) or (LA) it enables us to have a sufficiently large set of apartments where we know the triangle inequality is satisfied.  This in turn allows us to work (via Weyl-compatibility) to deduce the triangle inequality on the whole space.

\begin{remark}\label{Rem_history}
The original axiomatic definition of affine $\R$-buildings is due to Jaques Tits, who defined the ``syst{\`e}me d'appartements'' in \cite{TitsComo} by listing five axioms. The first four of these are precisely axioms (A1)--(A4) as presented above. His fifth axiom originally reads different from ours but was later replaced with what is now axiom (A5) as written in Definition~\ref{Def_LambdaBuilding}. One can show that in case of $\R$-buildings axioms (A6) follows from (A1)--(A5). 
In appendix 3 of  \cite{Ronan} it is stated that (A5) and (A6) may be used alternatively. The actual definition of an $\R$-building Mark Ronan gives is (\ref{a4}) of Theorem~\ref{MainThmB} with condition (TI) dropped.  It is however possible to construct examples of spaces satisfying (A1)--(A4) that vacuously satisfy (A6) but satisfy neither (A5) nor (TI). See Section~\ref{Sec_example} for details. Ronan \cite{Ronan} implicitly seems to assume condition (TI) in addition to (A1)--(A4) and (A6).  \newline
Anne Parreau shows \cite{Parreau} that in case of $\R$-buildings the conditions of Theorem~\ref{MainThmB} are also equivalent to (A1)--(A4) and (LA). It is not clear to us whether this may hold in the general setting. 
%\marginpar{boils down to prove that CO is equiv to A4 given LA} 
One might also ask whether for $\R$-buildings the weaker assumption (A4) and (aLA) might suffice. \newline
Guy Rousseau takes (A1),(A2) and (LA) as definition of $\R$-buildings\footnote{He calls them Euclidean buildings.} in \cite{Rousseau}. This definition seems to be equivalent for the underlying metric space but may allow for more general apartment systems. 
\end{remark}

\subsection{An application}
 
One can view Euclidean buildings as the subclass of affine $\Lambda$-buildings where $\Lambda=\R$ and where the translational part $T$ of the affine Weyl group equals the co-root-lattice spanned by a crystallographic root system, or is the full translation group of an apartment in the non-crystallographic case. 
Concerning the metric structure however the following difficulty arises when doing so. 
The metric used on the model space of a Euclidean buildings is usually the Euclidean one. Compare for example \cite{Parreau} or Kleiner and Leeb \cite{KleinerLeeb}. The natural metric on the model space of an affine $\Lambda$-building is however defined in terms of the defining root system $\RS$ (see \cite{Bennett, Diss}) and is a generalization of the length of translations in apartments of simplicial affine buildings. This length function on the set of translational elements of the affine Weyl group is defined with respect to the length of certain minimal galleries and differs from the Euclidean metric. 

The question arising is the following: Let us assume that $X$ is an affine building with metric $d$, which is induced by a metric $d_\MS$ on the model space. Let $d'_\MS$ be another metric on the model space which induces a second distance function $d'$ on $X$. 
Does $d'$ satisfy the triangle inequality? And is $(X,d')$ an affine building? In order to answer these questions one has to understand whether for $d'$ the retractions appearing in (A5) do exist and are distance diminishing. It turns out that the answer is ``yes'', so long as the metrics involved are Weyl compatible and topologically equivalent on $\MS$. 
 
\begin{thm}\label{thm_cor}
Let $(X,\App)$ be an affine $\Lambda$-building with  Weyl compatible metric $d$. Then every Weyl compatible metric $d'$ on $\MS$ that is topologically equivalent to $d$ on $\MS$ extends to a metric on $X$. In particular ``being a building'' only depends on the equivalence class (in the class of Weyl compatible metrics) of $d$, not on the metric itself.
 \end{thm}
\begin{proof}
For any Weyl compatible metric $d'$, that is topologically equivalent to $d$, axioms (A6) and (A1) to (A4) are still satisfied, since these axioms do not contain conditions on the metric. In addition, since the metrics $d$ and $d'$ are topologically equivalent on the model space, they define the same germs of Weyl chambers.  In particular, conditions (GG), (LA), and (aLA) are all satisfied for $d'$ if they are satisfied for $d$.  Since condition (CO) is independent of the metric, it follows that both statements (\ref{a8}) and (\ref{a9}) of \ref{MainThmB} are true for $d'$ also.  Consequently $(X,\App,d')$ is a $\Lambda$-building and the metric induced by $d'$ satisfies (TI).
\end{proof}

%%%%%%%%%%%%%%%%%%%
\section{Having a large atlas (LA)}\label{Sec_LA}

Suppose that $(X,\App)$ is a pair satisfying axioms (A1) to (A3). 
Recall that the germs of Weyl simplices based at a vertex $x$ are partially ordered by inclusion. A germ $\Delta_x S_1$ of a Weyl simplex $S_1$ is contained in the germ $\Delta_xS_2$ of a different Weyl simplex $S_2$ if there exist $x$-based representatives $S'_i$ of $\partial S_i, i=1,2$ contained in a common apartment such that $S_1'$ is a face of $S_2'$. The \emph{residue} $\Delta_xX$ of $X$ at $x$  is the set of all germs of Weyl simplices based at $x$.

\begin{thm}\label{Thm_residue}
If $(X,\App)$ satisfies (A1)--(A3) and property (GG), then $\Delta_xX$ is a spherical building of type $\RS$ for all $x$ in $X$. 
\end{thm}
\begin{proof}
We verify the axioms of the definition of a simplicial building, which can be found on  page 76 in \cite{Brown}.
It is easy to see that $\Delta_xX$ is a simplicial complex with the partial order defined above. It is a pure simplicial complex, since each germ of a face is contained in a germ of a Weyl chamber. The set of equivalence classes determined by a given apartment of $X$ containing $x$ is a sub-complex of $\Delta_xX$ which is a Coxeter complex of type $\RS$. Hence we define those to be the apartments of $\Delta_xX$. Therefore, by definition, each apartment is a Coxeter complex. 
Two apartments of $\Delta_xX$ are isomorphic via an isomorphism fixing the intersection of the corresponding apartments of $X$, hence fixing the intersection of the apartments of $\Delta_xX$ as well. Finally due to property (GG) any two chambers are contained in a common apartment and we can conclude that $\Delta_xX$ is a spherical building of type $\RS$.
\end{proof}

\begin{corollary}
Suppose $(X,\App)$ is an affine $\Lambda$-building. Then $\Delta_xX$ is independent of $\App$. 
\end{corollary}
\begin{proof}
Let $\App'$ be a different system of apartments of $X$. We will denote by $\Delta$ the spherical building of germs at $x$ with respect to $\App$ and  by $\Delta'$ the building at $x$ with respect to $\App'$. Since spherical buildings have a unique apartment system the buildings $\Delta$ and $\Delta'$ are equal if they contain the same chambers. Let $c\in\Delta'$ be a chamber; we will show $c\in\Delta$. Let $d$ be a chamber opposite $c$ in $\Delta'$ and $a'$ the unique apartment containing both. Then $a'$ corresponds to an apartment $A'$ of $X$ having a chart in $\App'$ and there exist $\App'$-Weyl chambers $S_c$, $S_d$ contained in $A'$ representing $c$ and $d$, respectively. Choose a point $y$ in the interior of $S_c$ and let $z$ be contained in the interior of $S_d$. By (A3) there exists a chart $f\in \App$ such that the image $A$ of $f$ contains $y$ and $z$. By \cite[Prop.~6.2]{Convexity2}  the point  $x$ is contained in $A$. And by construction the unique $x$-based Weyl chamber in $A$ 
which contains $y$ has germ $c$ and the unique $x$-based Weyl chamber in $A$ containing  $z$ has germ $d$. Thus  $c\in\Delta$.  Interchanging the roles of $\Delta$ and $\Delta'$ above we have that they contain the same chambers.   Hence $\Delta=\Delta'$.
\end{proof}

It is possible to weaken the assumptions of this corollary a bit. However in order to apply \cite[Proposition 6.2]{Convexity2} in its proof we need to assume at least (A1)--(A3), (A5) and (SC).

\begin{lemma}\label{Lm:secgm sec}
Assume in addition to (A1)--(A3) that $(X,\App)$ satisfies (GG) and (CO) or, alternatively,  that (A4), (TI),  and (SC) are satisfied. Let $S$ and $T$ be two $x$-based Weyl chambers. Then there exists an apartment containing $T$ and a germ of $S$ at $x$.  
\end{lemma}
\proof In case that (GG) and (CO) are satisfied, the proof is as in \cite[Prop 1.15]{Parreau}.  
Suppose now that (A4), (TI), and (SC) hold, and obtain from Theorem~\ref{thm:infty} that these axioms are enough to see that $\partial_\App X$ is a spherical building. We write $d(S,T)$ for the Weyl group-valued Weyl chamber-distance of $\partial S$ and $\partial T$ in the building $\partial_\App X$ at infinity. Taking $\ell$ to be the usual length function on Coxeter groups, we use  $\ell(d(S,T))$ to denote the length of a minimal gallery connecting the chambers $\partial S$ and $ \partial T$. 

Starting with an apartment $A_0$ containing $S$ we will now construct an apartment containing $T$ and a germ of $S$ at $x$.

By axiom (A4) there exists an apartment $A$ containing sub-Weyl chambers $S'$ of $S$ and $T'$ of $T$ with $\ell(d(S',T'))= \ell(d(S,T))$.  Replace $S'$ and $T'$ by Weyl chambers in $A$ based at a common vertex $x'\in S$, consider a minimal-length sequence 
$$
S'=S_0,\dots,S_n=T',
$$
of $x'$-based Weyl chambers.

By construction $A_0$ contains $S'$. Let $j\in\{0, 1,\ldots, n-1\}$ be minimal such that $A_0$ does not contain a sub-Weyl chamber of $S_{j+1}$.  If such a $j$ does not exist $T'$ has a sub-Weyl chamber $T''$ contained in $A_0$.  But then, as $x\in A_0$, by convexity (i.e. axiom (A2)) it follows that $T\subset A_0$ and there is nothing to prove.  

Suppose $S_j\subset A_0$ but $S_{j+1}$ has no sub-Weyl chamber contained in $A_0$ for some $j\leq n-1$.
Then by construction there exists a Weyl chamber $S_{j+1}'$ parallel to $S_{j+1}$ in $A$ such that $S_{j+1}'\cap A_0$ is a panel $P$.  Moreover $P$ is parallel to a panel of $S_j$ (and $S_{j+1}$).  Let $M$ be the hyperplane of the apartment $A_0$ that is spanned by $P$.
% 
% (The idea below is that
% the two apartments $A_0$ and $A_{j+1}$ have apartments at infinity, and
% what you follow is the chambers there, which satisfy the conditions of two
% apartments in a spherical building that agree on a half apartment).

From (SC) we deduce the existence of an apartment $A_{j+1}$ containing the Weyl chamber $S_{j+1}'$ and the germ $\Delta_xS$ (since the germ $\Delta_xS$ must lie in one of the half-apartments of $A_0$ determined by the hyperplane spanned by $P$) such that $A_{j+1}\cap A_0$ is a half-apartment.

There are two possible cases for $S$: 
(1) All of $S$ lies on the same side of the hyperplane $M$ as does the
germ $\Delta_x S$, in which case $S$ is contained in $A_{j+1}$. In this case we replace $A_0$ by $A_{j+1}$. 
(2) The interior of the Weyl chamber $S$ intersects the wall $M$ (i.e.\ $S$ does not, in particular, have a panel contained in that wall), in which case there is a sub-sector $S''$ of $S$ lying entirely on the opposite side of the wall $M$ from
$\Delta_x S$ inside $A_0$. 
We further denote by $\tilde S$ the unique Weyl chamber in $A_{j+1}$ with germ $\Delta_x\tilde S=\Delta_x S$. 
Since $A_{j+1}$ and $A_0$ agree on a half-apartment, there is a isomorphism $g$ from $A_0$ to $A_{j+1}$
preserving the intersection $A_0\cap A_{j+1}$, which necessarily preserves
the sector germ $\Delta_x S=\Delta_x \tilde{S} $.  Moreover, as $S_{j+1}$
is not in $A_0$, $S_{j+1}$ has a sub-Weyl chamber that lies on the opposite side of the wall $M=g(M)$
from $\Delta_x S$ in $A_{j+1}$.  But since $g$ is an isomorphism, it follows that
$g(S)=\tilde{S}$ (as $g(\Delta_x S )=\Delta_x\tilde{S}$.  Similarly,
since $g(P)=P$ and $\ell(d(S,S_j))<\ell(d(S,S_{j+1}))$, the Weyl chamber $S_j$ has a sub-Weyl chamber that lies on the
same side of the wall $M$ as $S''$ inside $A_0$.  Let $H$ denote the half apartment of $A_0$ defined by 
the wall $M$ that does not contain  $\Delta_x(S)$.  Then $g^{-1}(S_{j+1}')$ lies in $H$
and as it has panel $P$, it must intersect $S_j$ in a subsector.
Consequently, 
\begin{eqnarray*}
\ell(d(\tilde{S},S_{j+1}')) &=&
\ell(d(g^{-1}(\tilde{S}),g^{-1}(S_{j+1}'))) \\
   & = & \ell(d(S,S_j)) \\
   & < & \ell(d(S,S_{j+1})).
\end{eqnarray*}
Hence a minimal gallery connecting $\partial \tilde S$ and $\partial T$ is shorter than $n$.
Replace in case (2) the Weyl chamber $S$ by $\tilde S\subset A_{j+1}$ (which satisfies $\Delta_xS=\Delta_x\tilde S$) and the apartment $A_0$ by $A_{j+1}$. 

Inductively repeating the same argument with $\tilde S$ replacing $S$ we may finally find an apartment $A_n$ containing $\Delta_xS$ and $S_n=T$. If in each step $A_{j+1}$ contains $S$ then $A_n$ contains both $S$ and $T$.

The Weyl valued distance $\delta(\Delta_xS, \Delta_xT)$ is $d(\bar S, T)$ for some Weyl chamber $\bar S$ obtained inductively from $S$ by the above procedure. In particular $\ell(\delta(\Delta_xS, \Delta_xT))\leq \ell(d(S,T))$ by the above inequality.  
\qed

\begin{corollary} \label{Lm:opp sec}
Let $(X,\App)$ be a pair satisfying axioms (A1)--(A4), (TI), and (SC).
If $S$ and $T$ are Weyl chambers of $X$ based at $y$, and $\delta(\Delta_yS,\Delta_yT)$ is maximal, then $S$ and $T$ are contained in a common apartment.
\end{corollary}
\proof  Since $\delta(\Delta_yS,\Delta_yT)$ is maximal, the proof of Lemma~{\ref{Lm:secgm sec}} implies that $\ell(\delta(\Delta_yS,\Delta_yT))=\ell(d(S,T))$.  However, in this case the lemma implies the existence of an apartment $A$ containing $S$ and $T$. 
\qed

See also \cite[5.23]{Diss} for the above corollary and \cite[5.15]{Diss} for the following proposition. 

\begin{prop}\label{Prop:sec sec}
Assume that in addition to (A1)--(A3) the pair $(X,\App)$ satisfies either (GG) and (CO) or, alternatively, axioms (A4), (TI) and (SC). Then (LA) is satisfied as well. 
\end{prop}
\begin{proof}
Let $\Delta_xS$ and $\Delta_yT$ be two germs of Weyl chambers $S$ and $T$. We begin by showing that $\Delta_xS$ and $y$ are contained in a common apartment.  By (A2), there exists an apartment $A$ containing $x$ and $y$.  Let $C$ be a Weyl chamber of $A$ based at $x$ containing $y$.  By Lemma~\ref{Lm:secgm sec} there exists an apartment $A'$ of $X$ containing $\Delta_xS$ and $C$.  Take $S'$ to be the Weyl chamber of $A'$ based at $y$ containing $\Delta_xS$.  Again by Lemma~\ref{Lm:secgm sec}, there exists an apartment $A''$ containing $\Delta_yT$ and $S'$.  Since $\Delta_xS\subset S'$, it follows that $A''$ contains $\Delta_yT$ and $\Delta_xS$ as desired.
\end{proof}

Another implication worth stating is the following:
\begin{lemma}\label{Lem_sFC-LA}
Suppose $(X,\App)$ satisfies axioms (A1)--(A3). Then (sFC) implies (LA).  
\end{lemma}
\begin{proof}
Let $\mu$ and $\eta$ be two germs of Weyl chambers based at $z$ and $x$ respectively. Let $F$ be a Weyl chamber based at $x$ that has germ $\eta$ and let $A$ be an apartment containing $F$.
Pick $y$ in the interior of $F$, which implies $\eta \subseteq \seg_A(x,y)$. Then from (sFC) we may conclude that the segment 
$\seg_A(x,y)$ is covered by a finite number of Weyl chambers based at $z$ all
of which are contained in a common apartment with $\mu$. But then one of
these Weyl chambers, call it $S$, must contain $\eta$ and there exists an apartment $B$ containing $S$ and $\mu$. Hence (LA). 
\end{proof}

%%%%%%%%%%%%%%%
\section{Exchange axioms}
\label{Sec_ECSC} 

In this section, we prove equivalence of the sundial configuration (SC), the exchange condition (EC) and axiom (A6) given that (A1)--(A4) and (BI) are satisfied. Recall that if a pair $(X,\App)$ satisfies (A1)--(A5) then (TI) is satisfied and that then by~\ref{thm:infty} condition (BI) is satisfied. So Propositions~\ref{Pr:A6 A6'} and \ref{Pr:A6' A6''} hold true if (BI) is replaced by either (A5) or (TI) .  Given an apartment $A$  we let $\partial A$ denote the associated apartment at infinity.  For notational convenience we will also use lower case letters for apartments in the building at infinity when not using the $\partial$-notation.

\begin{prop} \label{Pr:A6 A6'}
Let $(X,\App)$ be a pair satisfying conditions~(A1)--(A4) and (BI) and suppose that the metric on $\MS$ is Weyl compatible, then condition~(A6) and the exchange condition~(EC) are equivalent.
\end{prop}
\proof
First assume $(X,\App)$ satisfies (A6). Suppose $A_1=f_1(\MS)$ and $A_2=f_2(\MS)$ are two apartments of $X$ with $A_1\cap A_2=F$ a half-apartment.  Then $\partial A_1$ and $\partial A_2$ are apartments of $\partial X$ that intersect in a half-apartment with bounding wall $h$.  By spherical building theory, it follows that there exists an apartment $a_3$ whose chambers are the chambers of $(\partial A_1\oplus \partial A_2) \cup\, h$, further there exists an apartment $A_3$ of $X$ with $\partial A_3 = a_3$.  

Since $\partial A_1\cap a_3$ is a half-apartment and $A_1\cap A_3$ is closed convex by (A2), it follows that $A_1\cap A_3$ is a half-apartment.  Similarly $A_2\cap A_3$ is a half-apartment.  Condition~(A6) now implies that $A_1\cap A_2\cap A_3$ contains some element $x\in X$.  Since $x\in F$ and $a_3$ contains the chambers of $\partial A_1$ that are not in $\partial A_2$, it follows that $A_3$ contains $A_1\setminus {F}$.  Similarly $A_2 \setminus {F}\subset A_3$.  By convexity the bounding wall $H$ of ${F}$ is contained in $A_3$.  But now the convexity of $A_3$ implies that $x  \in H$  as otherwise the wall parallel to $H$ through $x$ would not separate points of $A_1\cap A_3$ and $A_2\cap A_3$.  This implies that the exchange condition (EC) holds.

Now assume that (A1)--(A4), (BI) and (EC) are satisfied, and let $A_1$, $A_2$, and $A_3$ be apartments of $X$ such that any two intersect in a half-apartment.  By way of contradiction, suppose $A_1\cap A_2\cap A_3=\emptyset$.  Let $F_{i,j}=A_i\cap A_j$ for $i,j\in\{1,2,3\}$.  Since $F_{1,2}\cap F_{1,3}=\emptyset$, it follows that if $F$ is a half-apartment of $A_1$ with $F_{1,2}\cap F$ contained in the boundary $H$ of $F$, then $F_{1,3}\cap F$ is again a half-apartment.  Now the exchange condition (EC) implies that there exists an apartment $A_4$ such that $A_4=(A_1\oplus A_2) \cup H_{1,2}$, where $H_{1,2}$ stands for the bounding wall of $F_{1,2}$.  Note that $F_{1,3}\subseteq A_4$, so that $\partial A_4$ consists of the same Weyl chambers as $\partial A_3$.  However, the apartments of $X$ are in one-to-one correspondence with the apartments of $\partial X$.  Therefore, $A_3=A_4$. 
\qed

\begin{prop} \label{Pr:A6' A6''}
Suppose  $(X,\App)$ is a pair satisfying conditions~(A1)--(A4) and (BI) and suppose that the metric on $\MS$ is Weyl compatible, then (EC) and (SC) are equivalent. 
\end{prop}
\proof  Suppose $(X,\App)$ satisfies conditions~(A1)--(A4) and (EC).  Let $A_1$ be an apartment and $S$ a Weyl chamber such that $S\cap A_1$ is a panel of $S$.  Then $\partial A_1\cap \partial S$ is a panel in $\binfinity X$. Therefore  there is an apartment $a_2$ of $\binfinity X$ such that $\partial S\in a_2$, and $\partial A_1\cap a_2$ is a half-apartment.  Let $A_2$ be the apartment of $X$ with $\partial A_2 = a_2$.  Since $\partial A_1\cap a_2$ is a half-apartment the intersection $A_1\cap A_2$ is a half-apartment.  We now apply condition~(EC) to obtain (SC).

Conversely, suppose $(X,\App)$ satisfies conditions~(A1)--(A4) and~(SC), and let $A_1$ and $A_2$ be apartments of $X$ intersecting in a half-apartment.  Let $S$ be a Weyl chamber of $A_2$ such that $S\cap A_1$ is a panel $P$ of $S$, and let $M$ be the wall of $A_1$ containing $P$.  By (SC), there exists an apartment $A_3$ containing $M$ such that $A_1\cap A_3$ is a half-apartment and $A_2\cap A_3$ is a half-apartment containing $M\cup S$. By convexity, it follows that $A_3=(A_1\oplus A_2)\cup M$ as desired.
\qed

From the previous two propositions we deduce:
\begin{corollary}\label{cor_SC}
Every building $(X,\App)$ satisfies (EC) and (SC). 
\end{corollary}

The following proposition is used in the proof of Theorem~\ref{MainThmB} in order to show that item (\ref{a8}) implies item (\ref{a10}).

\begin{prop}\label{Pr:EC}
Let  $(X,\App)$ be a pair satisfying axioms (A1)--(A3) and (CO) and assume that the germs at each vertex form a spherical building\footnote{By Theorem~\ref{Thm_residue} this is true assuming (GG), for example.}.
Then the exchange condition (EC) is satisfied.
\end{prop}
\begin{proof}
Let $A$ and $B$ be apartments intersecting in an half-apartment $F$. Let $x$ be a point contained in the bounding wall $H$ of $F$. By assumption $\Delta_xX$ is a spherical building. Therefore the union of $\Delta_x(A\setminus F)$,  $\Delta_x(B\setminus F)$ and $\Delta_xH$ is an apartment in $\Delta_xX$, which we denote by $\Delta_xA'$. 

We choose two opposite germs $\mu$ and $\sigma$ at $x$ that are contained in $\Delta_x((A\setminus F)\cup H)$ and $\Delta_x((B\setminus F)\cup H)$, respectively, and that both have a panel germ contained in $H$. Let $T$ be the unique Weyl chamber in $A$ having germ $\mu$ and let $S$ be the unique Weyl chamber in $B$ with germ $\sigma$. By construction $S$ and $T$ are opposite and both have panels contained in $H$. Condition (CO) implies that $S$ and $T$ are contained in a common apartment $A''$. Since two opposite Weyl chambers contained in the same apartment determine this apartment uniquely we can conclude that  $\Delta_xA''=\Delta_xA'$. Also, by convexity, $A''$ needs to contain $H$. We conclude that $A''\cap ((A\oplus B)\cup H)$ contains $S$, $T$, $H$  and $\Delta_xA'$. From axiom (A2), that is convexity of intersection of apartments, we deduce  $A''\cap (B\setminus F) = B\setminus F$ and $A''\cap (A\setminus F) = A\setminus F$ which implies that $A''\cap ((A\oplus B)\cup H) = A''$.
\end{proof}

%%%%%%%%%%%%%
\section{Local structure}\label{Sec_localStructure}

Suppose $(X,\App)$ is a pair equipped with a Weyl compatible metric $d$ that satisfies (A1)--(A4), (BI) and (A6) as in Condition~(\ref{a11}) of Theorem~\ref{MainThmB}.  Recall that a germ $\mu$ of a Weyl chamber $S$ at $x$ is \emph{contained in a set $Y$} if there exists $\varepsilon\in\Lambda^{+}$ such that $S\cap B_\varepsilon(x)$ is contained in $Y$.

\begin{prop}\label{Prop_tec16}
Let $c$ be a chamber in $\binfinity X$ and $S$ an $x$-based Weyl chamber in $X$. Then there exists an apartment $A$ such that $\Delta_xS$ is contained in  $A$ and  such that $c$ is a chamber of $\partial A$.
\end{prop}

The proof of the proposition above is as in \cite[Prop.~1.8]{Parreau}. It uses axioms (A1)--(A4) and (A6).  Be warned (in case you want to look up Parreau's proof) that axiom (A5) in \cite{Parreau} is called (A6) in our list of axioms. As a direct consequence of~\ref{Prop_tec16} we obtain the following corollary.

\begin{corollary}\label{Cor_GG}
Any such pair $(X,\App)$ satisfying Condition~(\ref{a11}) of Theorem~\ref {MainThmB} has property (GG).
\end{corollary}
\begin{proof}
For a pair of germs $\mu$ and $\nu$ at $x$ pick Weyl chambers $S$ and $T$ both based at $x$ with germs $\mu$ respectively $\nu$. By Proposition~\ref{Prop_tec16} there exists then apartment $A$ of $X$ containing $S$ (but then also $\mu$) and $\nu$. 
\end{proof}

By the previous corollary the pair $(X,\App)$ satisfies the assertion of Theorem~\ref{Thm_residue}, i.e.\  the germs at a fixed vertex form a spherical building. 

\begin{prop}\label{Prop_LA}
Given $(X,\App)$ as above, then $(X,\App)$ has a large atlas, i.e., has property (LA).
\end{prop}
\begin{proof}
We need to prove that if $S$ and $T$ are Weyl chambers based at $x$ and $y$, respectively, then there exists an apartment containing a germ of $S$ at $x$ and a germ of $T$ at $y$.

By axiom (A3) there exists an apartment $A$ containing $x$ and $y$. We choose an $x$-based Weyl chamber $S_{xy}$ in $A$ that contains $y$ and denote by $S_{yx}$ the Weyl chamber based at $y$ such that $\partial S_{xy}$ and $\partial S_{yx}$ are opposite in $\partial A$. Then $x$ is contained in $S_{yx}$. If $\Delta_yT$ is not contained in $A$ apply Proposition~\ref{Prop_tec16} to obtain an apartment $A'$ containing a germ of $T$ at $y$ and containing  $\partial S_{yx}$ at infinity. But then $x$ is also contained in $A'$. 

Let us denote by $S'_{xy}$ the unique Weyl chamber contained in $A'$ having the same germ as $S_{xy}$ at $x$.
Without loss of generality we may assume that the germ $\Delta_yT$ is contained in $S'_{xy}$. Otherwise $y$ is contained in a face of $S'_{xy}$ and we can replace $S'_{xy}$ by an adjacent Weyl chamber in $A'$ satisfying this condition.
A second application of \ref{Prop_tec16} to $\partial S'_{xy}$ and the germ of $S$ at $x$ yields an apartment $A''$ containing $\Delta_xS$ and $S'_{xy}$ and therefore $\Delta_yT$.
\end{proof}

Propositions \ref{Prop_A6} to \ref{Prop_liftGallery} below  are due to Linus Kramer. A proof of \ref{Prop_A6} can be found in \cite[Prop.~5.20]{Diss}. 

\begin{prop}\label{Prop_A6}
Suppose $(X,\App)$ satisfies condition~(\ref{a11}) of Theorem~\ref{MainThmB}.  Let $A_i$ with $i=1,2,3$ be three apartments of $X$ pairwise intersecting in half-apartments. Then $A_1\cap A_2\cap A_3$ is either a half-apartment or a hyperplane.
\end{prop}

\begin{prop}\label{prop_sc}
Given $(X,\App)$ as above, then it satisfies (SC).
\end{prop}
\begin{proof}
Let $A$ be an apartment in $X$ and $c$ a chamber not contained in $\partial A$ but containing a panel of $\partial A$. Then $c$ is opposite two uniquely determined chambers $d_1$ and $d_2$ in $\partial A$. Since any pair of opposite chambers is contained in a common apartment,  there exist apartments $A_1$ and $A_2$ of $X$ such that $\partial A_i$ contains $d_i$ and $c$ with $i=1,2$.  The three apartments $\partial A_1,\partial A_2$ and $\partial A$ pairwise intersect in half-apartments. 
\end{proof}

Axiom (A6) together with the proposition above implies that the three apartments of the sundial configuration  intersect in a hyperplane. 

For any point $x\in X$ one can define a natural projection $\pi:\binfinity X \to \Delta_xX$ from the building at infinity to the residue at $x$ as follows.
Let $c$ be a chamber at infinity. Then there exists a unique Weyl chamber $S$ based at $x$ such that $\partial S = c$. Let 
$\pi (c) =\Delta_xS$.

\begin{prop}\label{Prop_liftGallery}
Suppose the pair $(X,\App)$ satisfies condition~(\ref{a11}) of Theorem~\ref{MainThmB}.  Let $(c_0, \ldots, c_k)$ be a minimal gallery in the building at infinity, $x$ a point in $X$ and for all $i$ denote by $S_i$ the $x$-based representative  of $c_i$. If $(\pi_x(c_0), \ldots, \pi_x(c_k))$ is a minimal gallery in $\Delta_xX$, then there exists an apartment containing $\bigcup_{i=0}^k S_i$.
\end{prop}

This follows by induction on $k$ and from repetitive use of Proposition~\ref{prop_sc}, compare \cite[Prop.~5.22]{Diss}.

\begin{corollary}\label{Cor_CO}
Suppose $(X,\App)$ satisfies condition~(\ref{a11}) of Theorem~\ref{MainThmB}.   Then $(X,\App)$ satisfies (CO).
\end{corollary}
\begin{proof}
Let $S$ and $T$ be Weyl simplices opposite at $x$. Choose a minimal gallery $(c_0,c_1,\ldots, c_n)$ from $c_0=\partial S$ to $c_n=\partial T$ and consider the representatives $S_i$ of $c_i$ based at $x$. With $S_0=S$ and $S_n=T$  Proposition~\ref{Prop_liftGallery} implies the assertion.
\end{proof}

\begin{corollary}\label{Cor_COTI}
If $(X,\App)$ satisfies (A1)-(A4), (TI) and (A6), then it satisfies (CO) and (GG). I.e. Condition~(\ref{a4}) implies both conditions~(\ref{a7}) and~(\ref{a8}) of Theorem~\ref{MainThmB}.
\end{corollary}
\begin{proof}
By Theorem~\ref{thm:infty}, $(X,\App)$ satisfies (BI), and now we apply Corollary~\ref{Cor_CO} and Corollary~\ref{Cor_GG}, 
\end{proof}

%%%%%%%%%%%%%%%%%%%%%
\section{Retractions based at germs}\label{Sec_retractions}

Let $(X,\App)$ be a pair satisfying (A1) and (A2) and assume $\App$ is an almost large atlas (aLA) (implying (A3)).  Further fix an apartment $A$ in $X$ with chart $f \in \App$.

\begin{definition}\label{Def_vertexRetraction}
\index{vertex retraction}
Let $\mu$ be a germ of a Weyl chamber in $A$ and $y$ a point in $X$, then, by (aLA),  there exists a chart $g\in\App$ such that $y$ and $\mu$ are contained in $g(\MS)$. 
By axiom (A2) there exists $w\in\WT$ such that $g\vert_{g^{-1}(f(\MS))}=(f\circ w)\vert_{g^{-1}(f(\MS))}$.
Hence we can define
$$
r_{A,\mu}(y) = (f\circ w\circ g^{-1} )(y).
$$
The map $r_{A,\mu}$ is called \emph{retraction onto $A$ centered at $\mu$}.
\end{definition}

\begin{prop}\label{Prop_r}
With $(X,\App), A, f $ and $\mu$ as above the retraction $r_{A,\mu}$ is well defined and the restriction of $r_{A,\mu}$ to any  apartment $A'$ containing $\mu$ is an isomorphism onto $A$.
\end{prop}
\begin{proof} 
By  (A2) the map $r_{A,\mu}$ is well defined.
Since each map $w\in W$ preserves distances in the model space $\MS$ we have that 
$$d(y,z)=d(r_{A, \mu}(y),r_{A, \mu}(z))$$ for all $y,z\in X$ such that $y$, $z$, and $\mu$ are contained in a common apartment. \end{proof}

We now verify finite covering properties that will allow us to prove under certain conditions that the defined retractions based at germs are distance diminishing.  

\begin{prop}\label{Prop_sFCa}
Assume  $(X,\App)$ is a pair satisfying (A1)--(A3) and properties (GG) and (CO).  Then (sFCa) and (aLA) are  satisfied. 
\end{prop}
\begin{proof}
We first prove (FCa).  Let $z\in X$ and $A\in\App$ be given.  We need to show there exists a finite cover of $A$ by Weyl Chambers based at $z$.
In the case where $z$ is contained in $A$ this is obvious. Hence we assume that $z$ is not contained in $A$. By (A3), for all $p\in A$ there exists an apartment $A'$ containing $z$ and $p$. 
Let $S_+\subset A'$ be a $p$-based Weyl chamber containing $z$.  We denote by $\sigma_+$ its germ at $p$. By condition (GG) and Theorem~\ref{Thm_residue} the link of $p$ is a spherical building. Hence there exists a $p$-based Weyl chamber $S_-$ in $A$ such that its germ $\sigma_-$ is opposite $\sigma_+$ at $p$. 
By property (CO) the Weyl chambers $S_-$ and $S_+$ are contained in a common apartment $A''$.
Let $T$ be the unique $z$-based translate of $S_-$ in $A''$. Since $z\in S_+$ and $\sigma_+$ and $\sigma_-$ are opposite we have that $S_-\subset T$. As there are only finitely many chambers in $\partial A$ there are only finitely many Weyl chambers based at $z$ with boundary in $\partial A$. Thus (FCa) follows. 

In order to see that (sFCa) holds true we argue as follows. Fix a germ $\mu$ at $z$. 
Let $I$ be the (finite) index set of the $z$-based Weyl chambers $S_i$ with equivalence class in $\partial A$. By (FCa) we may conclude that $A\subseteq \bigcup_{i\in I} S_i$, where each $S_i$ is based at $z$.
Fixing $i$, (GG) implies there exists an apartment ${A}_i$ containing $\mu$ and $\Delta_zS_i$. Let $S_i^{op}$ be a Weyl chamber in ${A}_i$ whose germ is opposite $\Delta_zS_i$ (such a Weyl chamber may be chosen as the link of $z$ is a spherical building). Then  (CO) implies that there is a unique apartment $A'_i$ containing the union of $S_i$ and $S_i^{op}$. By construction the induced apartment $\Delta_z A'_i=\Delta_z A_i$ as they both contain the same pair of opposite chambers.  Thus $\Delta_z A_i'$ contains $\mu$, $\Delta_z S_i^{op}$ and $\Delta_zS_i$. Since  $A\subseteq \bigcup_{i\in I} S_i$ and each $S_i$ is contained in a common apartment, $A_i'$ with $\mu$, we have shown (sFCa). As (sFCa) implies (aLA) we have the assertion.  
\end{proof}

Next we establish a local version of the sundial configuration.

\begin{property}\label{cond}
Recall that $(X,\App)$ satisfies condition~(\ref{a6}) of Theorem~\ref{MainThmB} if it satisfies conditions (A1)--(A4), (TI) and the sundial configuration (SC).  For the remainder of this section, we will just refer to this as condition~(\ref{a6}).
\end{property}

\begin{lemma}\label{Lm:sup exch}
With $(X,\App)$ as in condition~(\ref{a6}) let $A$ be an apartment of $X$ and $\Delta_xS$ a germ of a Weyl chamber such that $\Delta_xS\cap A$ is a panel-germ $\Delta_xP$. Then there exist apartments $A'$ and $A''$ such that $\Delta_xS \in A'\cap A''$ and $A\subset A'\cup A''$.
\end{lemma}
\proof
By (SC) it suffices to show that there is a Weyl chamber $S'$ of $X$ intersecting $A$ in a panel such that $\Delta_xS'=\Delta_xS$.  Let $T$ be a Weyl chamber of $A$ based at $x$ having a panel $P'$ containing the panel-germ $\Delta_xP$.  By Lemma~\ref{Lm:secgm sec} there exists an apartment $B$ containing $T$ and $\Delta_xS$.  Let $S'$ be the Weyl chamber of $B$ having the panel-germ $\Delta_xS$.  Then $S'$ has $P'$ as a panel.  Moreover, by convexity, if $S'\cap A\ne P'$, then $\Delta_xS=\Delta_x{S'}\subset A$ contrary to our hypothesis. Therefore $S'\cap A=P'$ and by (SC) there exists apartments $A'$ and $A''$ such that $S'\subset A'\cap A''$ and $A\subset A'\cup A''$.
\qed

This exchange condition allows us to work with germs based at a common point, much as in the simplicial buildings case one works with chambers in a spherical residue. 

The following proposition shows that condition (\ref{a6}) is enough to conclude (sFCa). Note that the assumptions made in \ref{Prop_sFCa} differ from the ones here. 

\begin{prop} \label{Pr:opps}
Suppose $(X,\App)$ satisfies condition~(\ref{a6}). Then (sFCa) is satisfied. In particular given an apartment $B$ of $X$ and $\Delta_xS$ a germ of a Weyl chamber, then for every point $y\in B$ there exists a Weyl chamber $T$ of $B$ such that 
\begin{enumerate}
\item the $x$-based Weyl chamber $T'$ parallel to $T$ contains $y$, and
\item there exists an apartment $A$ of $X$ containing $T$ and $\Delta_xS$.
\end{enumerate}
\end{prop}
\proof  By Proposition~\ref{Prop:sec sec}, for every Weyl chamber $T$ of $B$ based at $y$, there exists an apartment $B'$ of $X$ containing $\Delta_yT$ and $\Delta_xS$. Let $S'$ be a Weyl chamber of $B'$ based at $y$ containing $\Delta_xS$.  Choose $T$ such that $\ell(\delta(\Delta_yT,\Delta_yS'))$ is maximal. 

If $\Delta_yT$ and $\Delta_yS'$ are not opposite (that is $\delta(\Delta_yT,\Delta_yS')$ is not the longest element of $\overline{W}$) then let $\Delta_yP$ be a panel germ of $\Delta_yT$ such that the wall $M$ of $B'$ through $\Delta_yP$ does not separate $\Delta_yT$ and $S'$.  In the apartment $B$ there exists a Weyl chamber $R$ such that $\Delta_yR$ shares $\Delta_yP$ with $\Delta_yT$ and such that $\Delta_yR\neq \Delta_yT$. 

Lemma~\ref{Lm:sup exch} implies, that there exists an apartment $B''$ containing $S'$ and $\Delta_yR$. 
Hence, since $\Delta_yT$ and $S'$ lie on the same side of $M$, by convexity the apartment $B''$ also contains $\Delta_yT$.  In $B''$, we then have $\ell(\delta(\Delta_yR,\Delta_yS'))=\ell(\delta(\Delta_yT,\Delta_yS'))+1$, contradicting the choice of $T$. 

Hence we may assume that $\Delta_yT$ and $\Delta_yS'$ are opposite.  By Corollary~\ref{Lm:opp sec}, there exists an apartment $A$ of $X$ containing $S'$ and $T$.  But $\Delta_xS\subset S'$, so that $\Delta_xS\subset A$.  Moreover, since $A$ contains $T$, take $T'$ to be the Weyl chamber based at $x$ parallel to $T$ (in $A$).  Since $T$ and $S'$ are opposite Weyl chambers and $x\in T$, it follows that $y\in T'$. As there are only finitely many Weyl chambers in $B$ we get (sFCa), completing the proof of the proposition. 
\qed

\begin{corollary}\label{Lm:fin cover}
Suppose that either  $(X,\App)$ is as in condition~(\ref{a6}) or is a pair satisfying (A1), (A2), and (sFC).  Let $A$ be an apartment of $X$ and $\mu$ a germ of a Weyl chamber.  Then for each pair of points $x,y\in A$ there exist closed convex sets $C_1$, \dots, $C_n$ in $A$ such that 
\begin{enumerate} 
\item\label{fincover 1} $\seg_A(x,y)\subset C_1\cup\dots\cup C_n$ and
\item\label{fincover 2} for each $i$ there is an apartment containing $C_i$ and  $\mu$.
\end{enumerate}
Assume (sFCa) in place of condition (sFC) then %for all points $x,y\in A$ 
there exist closed convex $C_i\subset A$ which satisfy~(\ref{fincover 2}) and are such that  $A=C_1\cup\dots\cup C_n$.
\end{corollary}
\proof
In the first case Proposition~\ref{Pr:opps} implies (sFCa). Let thus $S_1, \ldots S_n$ be the Weyl chambers provided by the finite cover condition (sFCa) (respectively (sFC)), set $C_i = S_i\cap A$ and the corollary follows. 
\qed

\begin{lemma}\label{lem:claim}
Let $x,y$ be two points in the model space $\MS$ and suppose $p$ is a point in $\seg_\MS(x,y)$.  Then
\begin{enumerate}
 \item\label{1} $\seg_\MS(p,y)\subset \seg_\MS(x,y)$ and
 \item\label{2}  $p\in\seg_\MS(q,r)$ for all $q\in \seg_\MS(x,p)$ and all $r\in\seg_\MS(p,y)$. 
\end{enumerate}
\end{lemma}
\begin{proof}
The first claim is a direct consequence of the following computation
\begin{eqnarray*}
d(x,y) &\leq & d(x,q)+d(q,y)\\
  &\leq & d(x,p)+d(p,q)+d(q,y)\\
  &=& d(x,p) +d(p,y) =d(x,y),
\end{eqnarray*}
where $q$ is an arbitrary point in $\seg_\MS(p,y)$.

To show (\ref{2}) observe first that 
\begin{eqnarray*}
d(x,y) &\leq & d(x,q)+d(q,y) \\
  &\leq & d(x,q)+d(q,p)+d(p,y)\\
  &=& d(x,p)+d(p, y)=d(x,y),
\end{eqnarray*}
where we use for the last two equalities that $q\in \seg_\MS(x,p)$ and $p\in \seg_\MS(x,y)$. This implies that $p\in\seg(q,y)$. 
Similarly obtain $p\in\seg(r,x)$. 

Directly from (\ref{1}) we get that $r\in\seg_\MS(x,y)$.  Further we may apply (\ref{1}) to $r$ in $\seg_\MS(p, y)$ and use the fact that $p\in\seg_\MS(q,y)$ to conclude that $r\in\seg_\MS(q,y)$. Combining all this we may conclude
\begin{eqnarray*}
d(p,q)+d(p,r) &=& d(x,r)-d(x,p) + d(q,y)-d(p,y)\\
  &=& d(x,y)-d(r,y) + d(q,y)-d(x,p)-d(p,y)\\
  &=& d(q,y)-d(r,y) = d(q,r).
\end{eqnarray*}
Hence the lemma.
\end{proof}

\begin{lemma}\label{lem:subsegments}
Suppose $(X,\App)$ satisfies condition~(\ref{a6}) or is a pair satisfying (A1), (A2) and (sFC). Suppose further that $d$ is a Weyl compatible metric. Let $x,y$ be points in an apartment $A$ and suppose $\mathcal{C}$ is a collection of closed convex sets satisfying the conditions of Corollary~\ref{Lm:fin cover}.
There exists then %a subset $\{C_i\}_{i=1}^k\subset\mathcal{C}$ and 
a sequence of points $x=x_0,x_1,\ldots,x_n=y$ such that
\begin{enumerate}
 \item\label{one} $x_i\in\seg_A(x_{i-1}, x_{i+1})$  for all $1\leq i\leq n-1$  and 
 \item\label{two} there exists $C\in \mathcal{C}$ with $x_{i-1},x_i\in C$ for all $1\leq i\leq n$.
\end{enumerate}
\end{lemma}
\begin{proof}
Identify $A$ with the model space $\MS$ and denote by $\mathcal{H}$ the finite set of all hyperplanes that support panels of the sets in $\mathcal{C}$ and that separate $x$ and $y$. Put $k\define\vert\mathcal{H}\vert$.  If $k=0$ then $\seg_\MS(x,y)$ is contained in one of the sets $C\in\mathcal{C}$ and we can put $n=1$, $x=x_0$ and $y=x_1$ to be the desired sequence. So suppose without loss of generality that $k\geq 1$.  

We will recursively define a sequence of points $p_{\{s_i\}_i}$ in the segment of $x$ and $y$ that satisfies conditions \ref{one} and \ref{two}. We order these points lexicographically by their indexing sequences $\{s_i\}_{i=1}^{k+1}$, with $s_i\in\{0,1\}$. Put $x=p_{000...0}$ and $y=p_{100...0}$. Now by $k$-step recursion we will introduce an additional point in each segment of subsequent points in the sequence. We will do this in such a way that the points added in step $i$ all satisfy $s_{i+1}=1$ and $s_{j}=0$ for all $i+2\leq j\leq k+1$.

For the first step start with an arbitrary $H\in\mathcal{H}$. Weyl compatibility then gives us a point $p=p_{010...0}$ in $\seg_\MS(x,y)\cap H$ which allows us to consider the two segments $\seg_\MS(x,p)$ and $\seg_\MS(p,y)$ in place of $\seg_\MS(x,y)$. 

Now in step $i$ (where $1\leq i\leq k$) we deal with at most $2^{i-1}$ segments of pairs of subsequent points $p_{\{s_i\}_i} < p_{\{t_i\}_i}$. Note that for each such pair $s_j=t_j=0$ for all $j>i$. 
Fix such a pair $s=p_{\{s_i\}_i}$ and $t=p_{\{t_i\}_i}$. 

If there is no hyperplane in $\mathcal{H}$ separating $s$ and $t$ we do not add a point in this segment. 
If there exist hyperplanes in $\mathcal{H}$ separating $s$ and $t$ pick one of them and call it $H$. Weyl compatibility implies that $\seg_\MS(s,t)\cap H\neq \emptyset$. Let $q$ be a point in this intersection and put $r_j=s_j$ for all $j\neq i+1$, $r_{i+1}=1$ and $p_{\{r_i\}_i}=q$. In particular $s<q<t$ with respect to the lexicographical ordering of the indices.  

By construction the resulting sequence satisfies (\ref{two}). Condition (\ref{one}) follows from Lemma~\ref{lem:claim}.\ref{2}.
\end{proof}

We are now ready to prove that the retractions under consideration are distance diminishing.

\begin{prop}\label{Prop_A5}
Let $(X,\App)$ satisfy condition~(\ref{a6}) or be a pair satisfying axioms (A1), (A2), and (sFC) and suppose further that the metric on $\MS$ is Weyl compatible. Then for all apartments $A$ and germs $\mu$ of Weyl chambers contained in $A$ the retraction $r_{A,\mu}$ defined in \ref{Def_vertexRetraction} is distance non-increasing. In particular $(X,\App)$ satisfies (A5).
\end{prop}
\begin{proof}
Assuming $X$ is as in condition~(\ref{a6}) we may deduce from Proposition~\ref{Prop:sec sec} that (LA) is satisfied and therefore retractions based at germs are well defined. Hence using  Proposition~\ref{Prop_r} in both cases we can conclude: if $B$ is an apartment  containing $\mu$, $y$ and $z$ then $d(y,z)=d(r_{A,\mu}(y),r_{A,\mu}(z))$, so the result holds true.  Now, suppose $y$ and $z$ are arbitrary.  By (A2) there exists an apartment $B$ containing $y$ and $z$.  By Corollary~\ref{Lm:fin cover} there exist closed convex sets $X_1$, \dots, $X_n$ such that $\seg_B(y,z)\subset \bigcup_{i=1}^n X_i$ and each $X_i$ is contained in a common apartment with $\mu$.  Via Lemma~\ref{lem:subsegments}, Weyl compatibility condition~(\ref{cond:metric2}) implies that there is a sequence of points
$y=y_0,y_1,\dots,y_k=z$ such that $y_{i-1},y_i\in X_{j_i}$ for some $j_1$, \dots, $j_k$ and $y_i$ is in the segment of $y_{i-1}$ and $y_{i+1}$ for $i=1,\dots,k-1$.  Then
\begin{eqnarray*}
d(y,z)&=&\sum_{i=1}^k d(y_{i-1},y_i) \\
  & = & \sum_{i=1}^k d(r_{A,\mu}(y_{i-1}),r_{A,\mu}(y_i)) \\
  & \ge & d(r_{A,\mu}(y_0),r_{A,\mu}(y_k)) 
   =  d(r_{A,\mu}(y),r_{A,\mu}(z)),
\end{eqnarray*}
where we use the triangle inequality for $d$ restricted to $A$ in the next to the last step.  Thus, $r_{A,\mu}$ is distance diminishing and hence a retraction with the required properties of (A5).
\end{proof}

We now wish to show that conditions (TI) and (SC) can replace conditions (A5) and (A6) in the definition of an affine $\Lambda$-building. 

\begin{prop} \label{Th:main} 
Suppose $(X,\App)$ satisfies axioms (A1)--(A4) and $d$ is a Weyl compatible metric.  Then conditions (A5) and (A6) together are equivalent to (TI) and (SC) together.  In other words:  in Theorem~\ref{MainThmB}  item (\ref{a1}) is equivalent to  (\ref{a6}). 
\end{prop}
\proof
We have shown in Corollary~\ref{cor_SC} that (SC) is satisfied by an affine $\Lambda$-building. Combining~\ref{Pr:A6 A6'} and~\ref{Pr:A6' A6''} we deduce that (A6) and (SC) are equivalent given (A1)--(A5). It remains hence to see that axioms  (A1)--(A4), (TI) and (SC) imply the retraction condition (A5), which follows from Propositions~\ref{Prop_r} and \ref{Prop_A5}. This completes the proof.
\qed

\begin{prop}\label{Th:tentwo}
Suppose $(X,\App)$ satisfies (A1)-(A4), (sFC) and (EC) and that $d$ is a Weyl compatible metric.  Then $(X,\App)$ satisfies (A5). In other words: in Theorem~\ref{MainThmB} item~(\ref{a10}) implies item~(\ref{a2}).
\end{prop}
\proof
This follows immediately from Proposition~\ref{Prop_A5}.
\qed

%%%%%%%%%%%%%%%%%%%%%%%%%%%%%%%%%%%%%
\section{Verifying (A4)}\label{Sec_A4}

Suppose that $(X,\App)$ is a pair satisfying axioms (A1)--(A3) and properties (GG) and (CO) and recall that by Theorem~\ref{Thm_residue}  these assumptions are enough to conclude that the germs at a given vertex form a spherical building. 

\begin{prop}\label{Prop_A4}
The pair $(X,\App)$ satisfies (A4). 
\end{prop}
\begin{proof}
Let $S$ and $T$ be two Weyl chambers in $X$. We will show that by passing to sub-Weyl chambers $S'$ and $T'$ we can find an apartment containing both $S'$ and $T'$.

From the proof of \ref{Prop_sFCa} we may deduce %(since  intersection of apartments are convex) 
that for any $x\in T$ there exists an apartment $B$ containing $x$ and a sub-Weyl chamber $\tilde S$ of $S$. Let $S_x$ be the unique $x$-based translate $S_x$ of $\tilde S$ in $B$ and write $T_x$ for the unique $x$-based sub-Weyl chamber of $T$. 

We denote by $\delta(x)$ the length of a minimal gallery from $\Delta_xS$ to $\Delta_xT$ in the spherical building $\Delta_xX$. 
Since the number of possible values for $\delta(x)$ is finite we may without loss of generality assume  that $x\in T$ is chosen such that $\delta(x)$ is maximal. 

Now replace $S$ by $S_x$ and $T$ by $T_x$. % where $x$ is chosen such that $\delta(x)$ is maximal. 
By  Lemma~\ref{Lm:secgm sec} there exists an apartment $A$ containing $T$ and a germ of $S$ at $x$ and we denote by $S'$ the $x$-based Weyl chamber in $A$ which is opposite $S$ at $x$. Property (CO) implies that there is an apartment $A'$ containing $S$ and $S'$. By (A2) the intersection $A'\cap T$ is a convex subset of $T$. 
Let $z$ be an arbitrary point in this intersection. % (possibly $z=x$).
The unique $z$-based Weyl chambers $S_z$ and $S_z'$ parallel to $S$ and $S'$, respectively, are both contained in $A'$. 
By construction the length of a minimal gallery from $\Delta_zS_z$ to $\Delta_zT_z$ is not greater than $\delta(x)$. On the other hand, since $T$ and $S'$ are both contained in the apartment $A$, we can conclude
$$
\delta_z(T_z, S_z') = \delta_x(T, S') = d-\delta_x(S,T) = d-\delta(x)
$$
where $d$ is the diameter of an apartment of $\Delta_xX$, that is the diameter of the spherical Coxeter complex associated to the underlying root system $\RS$. The function $\delta_x$ assigns to two $x$-based Weyl chambers the length of a minimal gallery connecting their germs in $\Delta_xX$.

The germ $\Delta_zT_z$ lies on a minimal gallery connecting the opposite germs $\Delta_zS_z$ and $\Delta_zS_z'$. Such a  minimal gallery is contained in the unique apartment containing $\Delta_zS_z$ and $\Delta_zS_z'$, which is $\Delta_zA'$. Therefore $\Delta_zT_z$ is contained in $\Delta_zA'$ as well. 
This allows us to conclude that $A'\cap T$ contains a germ of $T_z$. One can observe that $A'\cap T$ is a convex subset of $T$ containing $x$ which is open relative to $T$. Hence the Weyl chamber $T$ is contained in $A'$. Thus (A4) follows. 
\end{proof}

\begin{thm}\label{eightten}
Suppose $(X,\App)$ satisfies (A1)-(A3), (GG) and (CO) and that $d$ is a Weyl compatible metric.  Then $(X,\App)$ satisfies (A4), (aLA), (sFC) and (EC).  In other words, in Theorem~\ref{MainThmB}, item~(\ref{a8}) implies item~(\ref{a10}).
\end{thm}
\proof
By Proposition~\ref{Prop_A4}, axiom (A4) holds.  In addition, by Theorem~\ref{Thm_residue}, for all $x\in X$, $\Delta_xX$ is a spherical building. Thus Proposition~\ref{Pr:EC} implies $(X,\App)$ satisfies (EC).  Proposition~\ref{Prop_sFCa} implies that $(X,\App)$ satisfies (sFCa) which in turn implies (sFC).  Thus item~(\ref{a8}) implies item~(\ref{a10}).
\qed

\section{The building at infinity}\label{Sec_infinity}

It was known already to the first author \cite{BennettThesis} that $\Lambda$-buildings (with respect to the standard $\Lambda$-metric) posses a spherical building at infinity. While his proof used the full power of axiom (A5) Parreau \cite{Parreau} later gave a proof for $\R$-buildings using (TI) only. We will give a proof here for pairs $(X,\App)$ with Weyl compatible metrics and show:  

\begin{thm}\label{thm:infty}
Let the model space $\MS=\MS(\RS,\Lambda)$ be equipped with a Weyl compatible metric $d$. If $(X, \App)$ is a pair modeled on $\MS$ satisfying axioms (A1)--(A4) and the triangle inequality (TI), then $X$ has a spherical building $\partial_\App X$ at infinity which is of type $R$. 
\end{thm}

Observe that translation invariance implies the following lemma.  

\begin{lemma}\label{lem:translates}
For all $v\in\MS$  there exists $\lambda>0$ such that for all $x\in\MS$ we have $d(x,x+v)\leq \lambda$.  
\end{lemma}

From this and Weyl compatibility condition \ref{cond:metric} (\ref{cond:metric1}) we may deduce: 

\begin{lemma}\label{lem:parallel}
Two Weyl simplices contained in a common apartment are translates of one another if and only if they are at bounded distance. 
\end{lemma}

\begin{definition}\label{Def:parallel}
We say that two Weyl chambers $S$ and $T$ in $X$ are \emph{parallel} if and only if they share a common sub-Weyl chamber. In this case write $S\sim T$. 
Two Weyl simplices $P$ and $Q$  are called \emph{parallel}, denoted by $P\sim Q$,  if and only if there exist a sequence of Weyl simplices $P=P_0,P_1,\dots,P_n=Q$ such that $P_i$ and $P_{i+1}$ are translates of each other in some apartment of $X$ for $i=0,\dots,n-1$.
\end{definition}

The parallel class of a Weyl simplex $P$ is denoted by $\partial P$ and we write $\partial X$ for the set of parallel classes of Weyl simplices in $X$.

\begin{lemma}
Parallelism is an equivalence relation. 
\end{lemma}
\begin{proof}
Symmetry and reflexivity are clear. Transitivity follows by concatenating sequences.
\end{proof}

Let $P$ be a sector-panel of $X$, lying in an apartment $f(\MS)$.  We define the type of $P$ to be the type of $f^{-1}(P)$ in $\MS$.  The type map is well defined by condition (A2).  Moreover, since translation in $\MS$ preserves the type of a sector panel, the type of a parallel class of sector-panels is also well-defined.

\begin{lemma} \label{lem:parbounddistance}
If $(X,\App)$ satisfies (A1)-(A4) and (TI), then Weyl simplices $P$ and $Q$ are parallel if and only if $P$ and $Q$ are at bounded distance.  Moreover, if $d(x,Q)\le \lambda$ for all $x\in P$, then for every sub-Weyl simplex $P'$ of $P$, there exists a sub-Weyl simplex $Q'$ of $Q$ such that for all $y\in Q'$ there exists $z \in P'$ with $d(z,y)\le \lambda$.
\end{lemma}
\begin{proof}
Suppose $P$ and $Q$ are parallel.  Then there exists a sequence $P=P_0,P_1,\dots,P_n=Q$ of Weyl simplices such that $P_i$ and $P_{i+1}$ are translates in an apartment $A_i$.  By Lemma~\ref{lem:parallel}, there exist  $\lambda_1,\dots,\lambda_n$ in $\Lambda$ such that for all $x\in P_{i-1}$ and $y\in P_i$, we have $d(x,y)\le \lambda_i$.  By the triangle inequality, it follows for $x\in P$ and $z\in Q$ that $d(x,z)\le\lambda_1 +\ldots +\lambda_n$, and consequently $P$ and $Q$ are at bounded distance in $X$.  Moreover, since in each apartment the condition on sub-Weyl simplices holds, the triangle inequality ensures that it holds in the concatenation.

Conversely, suppose $P$ and $Q$ are at bounded distance.  Let $S_1$ and $S_2$ be Weyl chambers containing $P$ and $Q$ as faces respectively.  By (A4), there exist sub-Weyl chambers $S_1'$ and $S_2'$ of $S_1$ and $S_2$ such that $S_1'$ and $S_2'$ are in a common apartment.  Moreover, there exists translates $P'$ and $Q'$ of $P$ and $Q$ in $S_1'$ and $S_2'$ respectively.  By the Lemma~\ref{lem:parallel}, $P$ and $P'$ and $Q$ and $Q'$ are at bounded distance, and hence by (TI) $P'$ and $Q'$ are at bounded distance.  Applying Lemma~\ref{lem:parallel} again, it follows that $P'$ and $Q'$ are translates of each other, implying that $P,P',Q',Q$ is a sequence of Weyl simplices showing that $P$ and $Q$ are parallel.
\end{proof}

We note that in the proof we also showed that in the case where $(X,\App)$ satisfies (A1)-(A4) and (TI), any two parallel Weyl simplices lie in a sequence  of length at most four of parallel simplices pairwise contained in apartments.  

\subsection*{Proof of Theorem~\ref{thm:infty}}

%%% CB 8/25 - is it (B1') in Brown that lets us look at Chamber and panel?  3.11 in Ronan, but he does not cite these as axioms, although Parreau seems to suggest that he does. %%%
We will prove that $\partial X$ is a building at infinity of $X$ by applying 3.11 in \cite{Ronan}. %showing that axioms (B0), (B1) and (B2') as stated in \cite{Brown}.

The set of parallel classes of Weyl simplices is ordered by inclusion and forms a chamber complex with the chambers being  equivalence classes of Weyl sectors. It is easy to see that each apartment in $X$ corresponds to a unique sub-complex of $\partial X$ isomorphic to a Coxeter complex of type $R$. Axiom (A4) immediately implies that any pair of simplices in $\partial X$ is contained in a common apartment. 

If a chamber $c\in\partial X$ is contained in two apartments $\partial A_1$ and $\partial A_2$ of $\partial X$  then there are two Weyl chambers $S_i\subset A_i$ both having the same boundary $c$, i.e.\ $c=\partial S_i$, $i=1,2$. Hence $S_1$ and $S_2$ have a sub-Weyl chamber $S$ in common. From (A2) we obtain that there is thus an isomorphism $f:A_1\to A_2$ fixing the intersection $A_1\cap A_2 \supset S_1\cap S_2$. The map $\partial  f$ induced by $f$ on the boundary is an isomorphism of Coxeter complexes fixing $\partial A_1 \cap \partial A_2$, we now use condition~\ref{cond:metric} (\ref{cond:metric1}).  Suppose $P_1\in A_1$ and $P_2\in A_2$ are parallel Weyl simplices (i.e., are such that $\partial P_1=\partial P_2$).  We must show $f(\partial P_1) = \partial P_2$.  Note that for all $x\in S_1\cap S_2$ and $y\in P_1$, $d(x,y)=d(f(x),f(y))=d(x,f(y))$.  
Moreover, by Lemma~{\ref{lem:parbounddistance}}, since $P_1$ and $P_2$ are parallel, there exists a $\lambda\in \Lambda$ 
such that for every sub-sector panel $P_1'\subset P_1$  there exists a sub-sector panel $P_2'\subset P_2$ such that for all $y\in P_2'$, there exists an $x\in P_1'$ such that $d(x,y)\le\lambda$.  By (TI), this implies that 
$$
|d(v,y)-d(v,f(x))| = |d(v,y)-d(v,x)| \le \lambda 
$$
for all $v\in S$.
Since $P_2$ and $f(P_1)$ and $S_1\cap S_2$ all lie in $A_2$ and are of the same type, condition~\ref{cond:metric} (\ref{cond:metric1}) implies that $f(P_1)$ and $P_2$ must be translates and therefore are parallel, so that $f(\partial P_1)=\partial(f(P_1))=\partial P_2$ as desired.
\qed

\section{Weyl compatibility}\label{Sec_WeylCompatibility}

By the \emph{standard $\Lambda$-metric} on the model space $\MS=\MS(\RS,\Lambda)$ we mean the $\Lambda$-valued metric as defined in \cite[Def. 4.16]{Diss} which is essentially the one introduced in the first author's thesis \cite{BennettThesis}.
For arbitrary $x$ and $y$ in $\MS$ we put
$$d(x,y)= \sum_{\alpha\in\RS^+} \vert\langle y-x, \alpha^\vee\rangle\vert  $$
where $\RS^+$ is the set of positive roots in the defining root system $\RS$.

Recall that the fundamental Weyl chamber $\Cf\subset \MS$, where  $\Cf=\{x\in\MS \;\vert\; \langle x,\alpha^\vee\rangle \geq 0 \text{ for all } \alpha\in \RS^+\}=\{x\in\MS \;\vert\; \langle x,\alpha^\vee\rangle \geq 0 \text{ for all } \alpha\in B\}$ with $B$ a basis of the defining root system $\RS$. Observe that each bounding wall of $\Cf$ is  ``perpendicular'' to an element of $B$. 
For points $x,y$ with $y-x  \in C_f$ one can show that 
$$
d(x,y) = 2\langle y-x , \rho^\vee \rangle,
$$ 
where $\rho^\vee=\frac{1}{2}\sum_{\alpha^\vee \in (R^\vee)^+} \alpha^\vee$ (see \cite[Section 1]{Diss}).

\begin{lemma}
The standard $\Lambda$-metric on $\MS$ is Weyl compatible.
\end{lemma}
\begin{proof}
Let us first verify condition~\ref{cond:metric} (\ref{cond:metric2}).
Let $x$ and $y$ be the two points separated by a wall $M$.  By $W$-invariance of the metric, without loss of generality, we may assume 
that $x$ is at the origin and $y\in \Cf$. Thus, the convex hull of $x,y$ is the set of all points such that the hyperplane coordinates $z^i$, where $\{\alpha_1,\dots,\alpha_n\}$ is the basis for $\RS$, satisfy  $0\le z^i\le y^i$ for $i=1,\dots,n$. For the definition of hyperplane coordinates see \cite{Bennett}. Let $\alpha\in\RS$ be the root such that $M$ is parallel to $H_{\alpha,0}$. Let $\alpha=\sum_{i} p_i \alpha_i$ (where $p_i\in \Q$ for all $i$).  
There exist $\lambda$ such that $M=\{z\in\MS\;\vert\; z^{\alpha}=\lambda \}$, where $z^\alpha = \sum_{i} p_i\cdot z^i$.
Now, consider the following ``path'' from $x$ to $y$.
$$
x=x_0,x_1,x_2,\ldots,x_n=y,
$$ 
where $x_j^i=0$ if $j>i$ and $x_j^i=y^i$ if $i\le j$.  If any $x_j$ lies
on $M$ we are done as each $x_j$ is in the segment from $x$ to $y$.  If
not, then there exists a $j$ such that $M$ separates $x_{j-1}$ from $x_j$.
 Let $\lambda_{j-1}=\sum_{i<j} p_i\cdot y^i$.  Then
$x_{j}^\alpha=\lambda_j$ so that
$$
\lambda_{j-1 }< \lambda < \lambda_j = \lambda_{j-1} + p_j\cdot y^j.
$$
It follows that $0 < \frac{1}{p_j} (\lambda-\lambda_{j-1}) < y^j$.  Since
$\Lambda$ is a $\Q$-module, it follows that $\gamma=\frac{1}{p_j}
(\lambda-\lambda_j) \in \Lambda$.  We now take $z$ to be the point in
$\Sigma$ given by the hyperplane coordinates $z^i=y^i$ for $i<j$, $z^j=\gamma$, and $z^{i}=0$ if $i>j$.
By construction $z\in\seg(x,y)\cap M$. 

% Without loss of generality we may assume that $x=0$ and $y$ is contained in the fundamental Weyl chamber $\Cf$. 
% Any hyperplane $H$ separating $x$ and $y$ is then of the form $H=H_{\alpha, k}$ for some positive root $\alpha$ and $0<k\in\Lambda$. Fix such a hyperplane $H_{\alpha, k}$ and let $\rho$ be half the sum of the positive roots.  By \cite[Prop. 29]{Bourbaki4-6} we know that $\rho\in\Cf$ and that $\langle\rho, \beta^\vee\rangle =1$ for all basis elements $\beta$. Hence $\langle \rho,\alpha^\vee\rangle>0$ for all positive roots $\alpha$.
% 
% Let $k'$ be such that  $y\in H_{\alpha, k'}$ and put 
% $$
% \lambda_{\alpha,k}\define \frac{2}{(\alpha, \alpha)\langle \rho, \alpha^\vee\rangle}(k'-k).
% $$ 
% Then 
% $$
% \frac{(\alpha,\alpha)}{2} \langle \lambda_{\alpha,k}\rho, \alpha^\vee\rangle = \lambda_{\alpha, k}\frac{(\alpha, \alpha)}{2}\langle\rho,\alpha^\vee\rangle = k'-k
% $$
% and $p=y- \lambda_{\alpha, k}\rho$ is contained in $H_{\alpha, k}$. 
% %It remains to show $p\in\seg(0,y)$. 
% %
% We further have that 
% \begin{eqnarray*}
% d(0,p) &=& 2\langle y-\lambda_{\alpha, k}\rho, \rho^\vee \rangle=2\langle y, \rho^\vee \rangle- 2\langle \lambda_{\alpha, k}\rho, \rho^\vee \rangle \;\;\text{ and }\\
% d(p,y) &=& d(0,y-p)=d(0,\lambda_{\alpha, k}\rho)=2\langle \lambda_{\alpha, k}\rho, \rho^\vee \rangle.
% \end{eqnarray*}
% Hence $p\in\seg_\MS(0,y)$ which implies Weyl compatibility condition (\ref{cond:metric2}).

% Proof part 1 Weyl compatibility

We will give an outline of the proof of Condition~\ref{cond:metric} (\ref{cond:metric1}).  Our outline follows the detailed proof of \cite[Section~2.6]{Bennett}.  Inside the model space $\MS$, we will show that if there exists a $\lambda$ such that for all sub-sector panels $P_1'\subset P_1$ there exist a sub-sector panel $P_2'\subset P_2$ with the property that for all $y\in P_2'$ there exists an $x\in P_1'$, such that for all $v\in \Cf$ we have the inequality
$$
|d(x,v)-d(y,v)|\le \lambda
$$
Then $P_1$ and $P_2$ are parallel, or $P_1$ and $P_2$ are parallel to Weyl-panels on a common Weyl chamber.  (Thus if they are the same type, they must be parallel.)

Suppose $P_1$, $P_2$, and $\lambda\in\Lambda$ are given as above.
Without loss of generality, we may assume that $P_1$ and $P_2$ are based at the origin. Hence there exist Weyl-chambers $S_1$ on $-P_1$ and $S_2$ on $-P_2$ such that $\Cf\subset P_i+S_i$.  Note that if $S_1$ is parallel to $S_2$ (and in fact equal since everything is based at the origin) then we have the result.  Thus, assume $S_1\ne S_2$, in which case there is a root $\alpha$ such that the wall $H_{\alpha,0}$ separates $S_1$ and $S_2$.

Let $v,q\in \Cf$ be given such that 
\begin{enumerate}
\item $q\in v+S_1$ (and note that $v+S_1\subset P_1+S_1$),
\item $|(\alpha,q-v)|\ge 3\lambda$ for all $\alpha\in \RS$.
\end{enumerate}
(That such $v$ and $q$ exist is shown in 2.19 of \cite{Bennett}).)

Let $P_1'$ be a sub-Weyl simplex of $P_1$ contained in $v-S_1$, and let $P_2'$ be the associated sub-Weyl simplex guaranteed by our hypothesis.  Then, let $P_2''\subset P_2'\cap (q-S_2)$, and let $y\in P_2''$.  From the triangle inequality, we have $d(y,q)\le d(y,v)+d(q,v)$.  Since there exists at least one $\alpha$ that separates $S_1$ and $S_2$, there is at least one $\alpha$ for which $q$ and $y$ lie on the same side of the wall $M_{\alpha,v}$ through $v$ parallel to $M_\alpha$.   Thus,
$$
d(y,q) \le  d(y,v)+d(q,v)-3\lambda.
$$
Let $N\ge 3\lambda$ be such that $d(y,q)=d(y,v+d(q,v))-N$.

Let $x\in P_1'$ be the point such that $|d(y,v)-d(x,v)|\le \lambda$.  Since $x\in P_1'\subset v-S_1'\subset q-S_1'$, we have
$d(x,q)=d(x,v)+d(v,q)$.  

Putting these together we obtain
\begin{eqnarray*}
|d(y,q)-d(x,q)| &=& \left| (d(y,v)+d(q,v)-N) - (d(x,v) + d(v,q))\right| \\
 & = & \left|d(y,v)-d(x,v)-N \right| \\
& \ge & N - |d(y,v)-d(x,v)| \\
& \ge & 3\lambda - \lambda >\lambda.
\end{eqnarray*}

This contradicts the hypothesis, hence $P_1$ and $P_2$ must either be parallel or lie on a common Weyl-chamber.   Moreover, if they are of the same type, then they must be parallel.

\end{proof}

The proof for the Euclidean metric being Weyl compatible in the case $\Lambda=\R$ is more straightforward.
\begin{lemma}
If $\Lambda=\R$, the Euclidean metric is Weyl compatible.
\end{lemma}
\begin{proof}
We begin by noting that condition~\ref{cond:metric} (\ref{cond:metric2}) follows from Hilbert's axioms on $\R$, given that a wall separates $\R^n$.  For condition~\ref{cond:metric} (\ref{cond:metric1}), we note that given any ray $r$ lying in $\Cf$ and any $\lambda\in\R^+$, the statement that $|d(v,y)-d(v,x)|\le \lambda$ implies that $y$ lies between two planes perpendicular to $r$, each lying distance $\lambda$ from $x$ (using the limit as $v\in r$ goes to infinity).  Since $\Cf$ contains a basis for $\MS$, it follows that $y$ lies in a parallelepiped about $x$, with the distance between any opposite faces being $2r$.  This implies that $y$ is at bounded distance from $x$.  Thus if $P_1$ and $P_2$ satisfy the hypothesis of condition~\ref{cond:metric} (\ref{cond:metric1}), they must be parallel (and in fact, we did not need that they were of the same type).  
\end{proof}

Given the much easier proof in the Euclidean metric case, one might wonder whether we need that $P_1$ and $P_2$ are of the same type in general.  Unfortunately, for the standard $\Lambda$-metric in the case of $\RS$ being the root system associated to the $A_2$ Coxeter diagram, the sector-panels opposite the fundamental sector panels of $\Cf$ satisfy condition~\ref{cond:metric}(\ref{cond:metric1}) except for the fact that they are of different types.

We believe that it might be possible to replace condition~\ref{cond:metric} (\ref{cond:metric2}) with a more natural concept of a geodesic metric.  Recall that a $\Lambda$-metric $d$ on a space $Y$ is \emph{geodesic} if for all $x,y\in Y$ there exists an isometric embedding of the interval $[0,d(x,y)]\subset\Lambda$ into $(Y, d)$. It is clear that any geodesic $\R$-metric on $\MS$ is Weyl compatible.

The standard $\Lambda$-metric is also geodesic. 

\begin{lemma}\label{Lem:geodesic}
The standard $\Lambda$-metric on $\MS$ is geodesic.  
\end{lemma}
\begin{proof}
Let $d$ denote the standard $\Lambda$-metric and suppose without loss of generality that $x=0$ and $y$ is contained in the fundamental Weyl chamber $\Cf\subset \MS$, where  $\Cf=\{x\in\MS \;\vert\; \langle x,\alpha^\vee\rangle \geq 0 \text{ for all } \alpha\in \RS^+\}=\{x\in\MS \;\vert\; \langle x,\alpha^\vee\rangle \geq 0 \text{ for all } \alpha\in B\}$ with $B$ a basis of the defining root system $\RS$. Observe that each bounding wall of $\Cf$ is  ``perpendicular'' to an element of $B$. 

For points $x,y$ with $y-x  \in C_f$ one can show that $d(x,y) = 2\langle y-x , \rho^\vee \rangle$, where $\rho^\vee=\frac{1}{2}\sum_{\alpha^\vee \in (\RS^\vee)^+} \alpha^\vee$, compare \cite[Section 1]{Diss}.  So in our case $d(0,y)=2\langle y, \rho^\vee \rangle$.  

Since the convex hull of a pair of points in $\MS$ equals their segment, see \cite[Prop. 2.13]{BennettThesis} we have that $\seg_\MS(0,y)\subset \Cf$.  

We write $\mathrm{proj}_{H_{\alpha_i, 0}}(p)$ for the projection $m_\alpha$ of a point $p$ onto the hyperplane ${H_{\alpha_i, 0}}$ perpendicular to $\alpha$. Suppose that the elements of the basis $B$ of $\RS$ are enumerated such that $B=\{\alpha_i\}_{i=1}^k$. We may define a sequence of projections (by backward induction) from $k+1$ to $1$ as follows: 
Put $y_{k+1}= y$ and define  
$$
y_{i}= \mathrm{proj}_{H_{\alpha_i, 0}}(y_{i+1}) \text{ for all } i=k, \ldots, 1. 
$$
With $y_0=0$ let $\gamma$ be the concatenation of the geodesic pieces $\gamma_i:y_{i} \rightsquigarrow y_{i+1}$, where $\gamma_i$ is defined as follows:
$$
\gamma_i(\lambda)=y_i + \lambda\alpha_i \text{ for all } \lambda\in [0, y_i^{\alpha_i}]  \text{ and all } i=0, \ldots, k.
$$
Here $x^\alpha=\frac{1}{2}\langle x, \alpha^\vee\rangle$ for all $x\in\MS$ and $\alpha\in \RS$, compare \cite[Prop. 4.12]{Diss}.
It is not hard to verify that indeed $d(0, y)=\sum_{i=0}^k d(y_i, y_{i+1})$ and that $\gamma: 0\rightsquigarrow y$ is a geodesic. 
\end{proof}

\section{Appendix: Sharpness of axiom (A5)}\label{Sec_example}

\begin{center}
by Koen Struyve
\end{center}

This section is devoted to the following question: Let $(X,\App)$ be a space modeled on $\A = \A (\RS,\Lambda)$, which satisfies axioms (A1)--(A4) and (A6). Is $(X,\cA)$ again an affine $\Lambda$-building?

We will construct examples of such spaces which vacuously satisfy (A6) but do not satisfy either (A5) or (TI), answering the above question in the negative.

\subsection{Some additional definitions}
The images in the model space $\A$ of the fundamental Weyl chamber $\cC_f$ under the spherical Weyl group $\overline{W}$ will be called \emph{vector Weyl chambers}. The unique fixed point in $\A$ of this group (so the basepoint of all the vector Weyl chambers) is denoted by $o$. The image of a vector Weyl chamber under a chart will again be called such.

We call the closed ball $\{x \in \A \vert d(o,x) \leq \lambda \}$ and their images under charts  \emph{centered balls with radius $\lambda$}.

Recall that if two Weyl chambers contain a common Weyl chamber we call them \emph{parallel}. If an atlas satisfies axiom (A2), then this relation forms an equivalence relation.

A space $(X,\cA)$ modeled on $\A = \A (\RS,\Lambda)$ will be called \emph{$\lambda$-admissible} if the following conditions are satisfied.
\begin{itemize}
\item[(T0)] No two injections in $\cA$ have the same image.
\item[(T1)] $(X,\cA)$ satisfies the axiom (A2).
\item[(T2)] If two different apartments intersect, then they either intersect in a single point contained in the centered balls with radius $\lambda$ of both apartments, or they intersect in a Weyl chamber contained in the interior of vector Weyl chambers in both apartments.
\item[(T3)] All the vector Weyl chambers in one parallelism class contain a common sub-Weyl chamber.%, not intersecting any centered ball of radius $l$.
\end{itemize}

Fix a sequence $(\lambda_i)_{i=1,2,\ldots}$ with $\lambda_i \in \Lambda$ such that $ 0< \lambda_1 < \lambda_2 < \dots$ and the sequence converges to infinity.

\subsection{Extension procedure}

In this section, we will construct from a given $\lambda_i$-admissible space $(X_i,\cA_i)$ (for $i\in \mathbb{N} \backslash \{0\}$) a $\lambda_{i+1}$-admissible space $(X_{i+1},\cA_{i+1})$, extending the previous one. %Iterating this two-step algorithm will produce a sequence of atlases such that the direct limit satisfies (A2), (A3) and (A4).

\subsubsection{Step 1: Covering pairs of points}

Let $P$ be the set of pairs of points (up to order) in centered balls with radius $\lambda_{i+1}$ in $X_i$ not yet covered by a common apartment. For each pair $p := (x,y)$ in $P$ we choose distinct points $x_p$ and $y_p$ in the centered ball of radius $\lambda_{i+1}$ of a copy $\A_p$ of the model space $\A$. Let $\pi_p$ be the canonical isometry from $\A$ to $\A_p$.

We now define $X'_i$ to be the union of the sets $X_i$ and $\A_p \setminus \{x_p,y_p\}$ for each $p:= (x,y)$ in $P$. The set of charts $\cA'_i$ is defined as the set of charts $\cA_i$ together with a chart 
$$
f_p: a \in \A \to \left\{\begin{array}{cl} \pi_p(a) & \mbox{if } x_p \neq \pi_p(a) \neq y_p  \\ 
x & \mbox{if } \pi(a_p) = x_p  \\ 
y & \mbox{if } \pi(a_p) = y_p  
\end{array}\right.  $$
for each $p:= (x,y)$ in $P$.

It is straightforward to verify that the newly obtained space $(X_i',\cA_i')$ satisfies conditions (T0) up to (T3) (with $\lambda =\lambda_{i+1}$), and hence is an $\lambda_{i+1}$-admissible space.

\subsubsection{Step 2: Covering pairs of sectors}
We will now extend $(X_i',\cA_i')$ to a $\lambda_{i+1}$-admissible space $(X_{i+1},\cA_{i+1})$.

Let $Q$ be the set of parallelism classes of Weyl chambers in $X_i'$. For each class $q\in Q$, consider all vector Weyl chambers in this parallelism class. We know that these vector Weyl chambers contain a common sub-Weyl chamber by condition (T3). Condition (T2) together with the definition of vector Weyl chambers then implies that there exists a sub-Weyl chamber which has no points in common with the centered balls of radius $\lambda_{i+1}$ in any apartment. Fix such a sub-Weyl chamber $S_q$. 

Let $R$ be the set of pairs (up to order) of parallelism classes of sectors not yet covered by a common apartment (meaning that there are no two elements, one of each class, contained in a common apartment).  For each such pair $r := (q_1,q_2)$ we pick two disjoint sub-Weyl chambers $S^1_r$ and $S^2_r$, both contained in vector Weyl chambers of a copy $\A_r$ of the model space $\A$, and such that $S^1_r$ and $S^2_r$ are disjoint with the centered ball of radius $\lambda_{i+1}$ in $\A_r$. Let  $\pi_r$ be the canonical isometry from $\A$ to $\A_r$. Let $\pi^1_r$ and $\pi^2_r$ be the canonical isometries from $S^1_r$ and $S^2_r$ to respectively $S_{q_2}$ and $S_{q_1}$. 

The set of points $X_{i+1}$ of the space we want to construct is the union of the sets $X_i'$ and $\A_r \setminus (S^1_r \cup S^2_r)$ for each $r$ in $R$. The set of charts $\cA_{i+1}$ is $\cA'_i$ extended with a chart

$$
g_r: a \in \A \to \left\{\begin{array}{cl} \pi_r(a) & \mbox{if } \pi_r(a) \notin S^1_r \cup S^2_r  \\ 
\pi^1_r(a) & \mbox{if }  \pi_r(a) \in S^1_r  \\ 
\pi^2_r(a) & \mbox{if }  \pi_r(a) \in S^2_r  \\ 
\end{array}\right.  
$$
for each $r$ in $R$.

We now claim that $(X_{i+1},\cA_{i+1})$ is $\lambda_{i+1}$-admissible space. Conditions (T0) and (T3) are automatically satisfied. In order to see conditions (T1) and (T2) note that for  any two points $a\in S_q$ and $b \in S_{q'}$, where $q,q'$ with $(q,q')\in R$\footnote{$q$ and $q'$ are not in the same apartment} are two distinct parallelism classes in $Q$, one has that $a$ and $b$ are not contained in a common apartment by condition (T2) for $(X_i',\cA_i')$.

\subsection{Direct limit and conclusion}

By repeating the extension procedure laid out in the previous sub-section one obtains sets of points $X_i \subset X_{i+1} \subset X_{i+2} \subset \dots$ and sets $\cA_i \subset \cA_{i+1} \subset \cA_{i+2} \subset \dots$ of injections. Let $$X_\infty = \bigcup_{j=i}^\infty X_{j} \text{ and } \cA_\infty = \bigcup_{j=i}^\infty \cA_{j} .$$

This direct limit yields a space $(X_\infty,\cA_\infty)$ modeled on $\A = \A (\RS,\Lambda)$. In order to satisfy axiom (A1) we replace $\cA_\infty$ by the set $\cA_\infty' = \{f \circ w \vert f\in \cF, w \in W\}$. The space $(X_\infty,\cA_\infty')$ satisfies axiom (A2) by condition (T1) for the intermediate steps. The repetition of the first and second step of the procedure implies that $(X_\infty,\cA_\infty)$ satisfies axioms (A3) and (A4) as well.

If the dimension of $\A$ is at least 2 then no two apartments intersect in a half-apartment by condition (T2) for the intermediary steps, so axiom (A6) is satisfied vacuously. So in this case we obtain a space which satisfies axioms (A1)--(A4) and (A6). However it cannot consists of more than a single apartment and satisfy axiom (A5) at the same time because if it would, then there would be apartments intersecting in half-apartments (see for example~\cite[Prop.~1.7]{Parreau}). By Theorem~\ref{MainThmB} it cannot satisfy (TI) either then.

A more direct way to obtain an example which does not satisfy (TI) is to start from a $\lambda_1$-admissible space (for a suitable choice of $\lambda_1$) which does not satisfy (TI). An example would be three apartments glued pairwise together at a (distinct) point. The three gluing points form a triangle, for which the side lengths can be chosen such that they violate the triangle inequality.

\phantomsection
\renewcommand{\refname}{Bibliography}
\bibliography{literaturliste-axioms}
\bibliographystyle{alpha}

\end{document}